\documentclass[10pt,reqno]{amsart}
\usepackage{amsmath,amstext,amsfonts,%
enumitem,%
extpfeil,graphicx,mathrsfs,%
mathtools,setspace,tikz,xparse}
\usepackage[alphabetic]{amsrefs}
\usepackage[unicode]{hyperref}
\usepackage[utf8]{inputenc}
\usepackage[T1]{fontenc}
\usepackage[normalem]{ulem}
\usepackage[cmtip,all]{xy}
\usepackage[light]{sourcesanspro}
\usepackage{lmodern}
\hypersetup{colorlinks=false,%
linkbordercolor={0.93 0.71 0.13},%
citebordercolor={1 0 .5},%
urlbordercolor={0.27 0.39 0.68}}%
\usetikzlibrary{calc,matrix}

\usepackage{filecontents}

\DeclareMathOperator{\End}{End}

\DeclareMathOperator{\Supp}{supp}
\DeclareMathOperator{\GL}{GL}
\DeclareMathOperator{\Sym}{Sym}
\DeclareMathOperator{\Img}{im}

\DeclareMathOperator{\ind}{ind}

\newcommand{\m}[1]{\(\smash{#1}\)}

\newcommand{\Algn}[1]{\begin{align}\textstyle #1\end{align}}
\newcommand{\bb}[1]{\mathbb{#1}}
\newcommand{\QQ}{\bb{Q}}
\newcommand{\FF}{\bb{F}}
\newcommand{\ZZ}{\bb{Z}}

\newcommand{\FFp}{\overline{\bb{F}}_p}
\newcommand{\QQp}{\smash{\overline{\bb{Q}}_p}}
\newcommand{\QQptimes}{\smash{\overline{\bb{Q}}_p^\times}}

\newcommand{\ZZp}{\smash{\overline{\bb{Z}}_p}}

\newcommand{\defeq}{\vcentcolon=}
\newcommand{\SCR}{\mathscr}
\newcommand{\IE}{i.e.\ }

\newcommand{\HECKET}{\mathbf T}

\newcommand{\MU}{\mu}

\newcommand{\SIGMA}{\Sigma}

\newcommand{\LAMBDA}{\lambda}

\newcommand{\LESS}[1]{%
{\uline{#1}}}
\newcommand{\MORE}[1]{%
{\overline{#1}}}
\newcommand{\BIGO}[1]{\smash{\mathsf{O}\!\left({#1}\right)}}

\newcommand{\SIZEEE}{4.25pt}
\newcommand{\VERIFYIF}[1]{[{#1}]}

\newcommand{\TEMPORARYA}{\xi}

\newcommand{\TEMPORARYD}{\vartheta}%

\newcommand{\TEMPORARYI}{\nu}

\newcommand{\absdet}{|\!\det\!|}

\newcommand{\PADICVALUATION}{v_p}

\newcommand{\TEMPORARYF}{L}
\newcommand{\TEMPORARYDISC}{D}

\newcommand{\TILDE}{\widetilde}

\newcommand{\GEQSLANT}{\geqslant}
\newcommand{\LEQSLANT}{\leqslant}

\newcommand{\EQUIVZEROIDEAL}{\equiv_{\mathfrak I}}

\newcommand{\maxideal}{\mathfrak m}
\newcommand{\repnn}{\pi}

\newcommand{\SUBMODULE}{\mathrm{sub}}
\newcommand{\QUOTIENTI}{\mathrm{quot}}

\newcommand{\FILT}{\widehat}

\newcommand{\DELTAA}{\Delta}

\newcommand{\sqbrackets}[4]{{#1}\bullet_{#2,#3}{#4}}

\newcommand{\Hypergeometric}[5]%
{\smash{\null_{#1}\mathrm{F}_{#2}\big(%
\smash{\genfrac{}{}{0pt}{2}{#3}{#4}};{#5}\big)}}

\newtheorem{theorem}{Theorem}
\newtheorem{lemma}[theorem]{Lemma}

\newtheorem{corollary}[theorem]{Corollary}
\newtheorem{conjecture}{Conjecture}

\newtheorem{innercustomthm}{Theorem}
\newenvironment{customthm}[1]
  {\renewcommand\theinnercustomthm{#1}\innercustomthm}
  {\endinnercustomthm}
  \newtheorem{innercustomcor}{Corollary}
\newenvironment{customcor}[1]
  {\renewcommand\theinnercustomcor{#1}\innercustomcor}
  {\endinnercustomcor}

\makeatletter
\newcommand{\dashover}[2][\mathop]{#1{\mathpalette\df@over{{\dashfill}{#2}}}}
\newcommand{\fillover}[2][\mathop]{#1{\mathpalette\df@over{{\solidfill}{#2}}}}
\newcommand{\df@over}[2]{\df@@over#1#2}
\newcommand\df@@over[3]{%
  \vbox{
    \offinterlineskip
    \ialign{##\cr
      #2{#1}\cr
      \noalign{\kern1pt}
      $\m@th#1#3$\cr
    }
  }%
}
\newcommand{\dashfill}[1]{%
  \kern-.5pt
  \xleaders\hbox{\kern.5pt\vrule height.4pt width \dash@width{#1}\kern.5pt}\hfill
  \kern-.5pt
}
\newcommand{\dash@width}[1]{%
  \ifx#1\displaystyle
    2pt
  \else
    \ifx#1\textstyle
      1.5pt
    \else
      \ifx#1\scriptstyle
        1.25pt
      \else
        \ifx#1\scriptscriptstyle
          1pt
        \fi
      \fi
    \fi
  \fi
}
\newcommand{\solidfill}[1]{\leaders\hrule\hfill}
\makeatother

\newenvironment{Smatrix}%
    {\tiny (\begin{smallmatrix}}%
    {\end{smallmatrix})}
\newenvironment{Lmatrix}%
    {\footnotesize (\begin{smallmatrix}}%
    {\end{smallmatrix})}
\newcommand{\BLACKSQUARE}{\hspace{\stretch{1}}\m{\mathbin{\vcenter{%
    \hbox{\rule{.9ex}{.9ex}}}}}}
\renewenvironment{proof}[1]%
    {\emph{{#1}}}%
    {\BLACKSQUARE}

\makeatletter
\newcommand{\xhookdoubleheadrightarrow}[2][]{%
  \lhook\joinrel
  \ext@arrow 0359\rightarrowfill@ {#1}{#2}%
  \mathrel{\mspace{-15mu}}\rightarrow}
\renewcommand{\xtwoheadrightarrow}[2][]{%
  \ext@arrow 0359\rightarrowfill@ {#1}{#2}%
  \mathrel{\mspace{-15mu}}\rightarrow}
\def\@settitle{\begin{center}%
  \baselineskip14\p@\relax\normalfont\LARGE
    \uppercasenonmath\@title
  \@title\end{center}%
}
\def\do#1{\@ifdefinable#1{\newdimen#1}}
\do\@TempWidthOne
\do\@TempWidthTwo
\newcommand*\WidthInMathUnitsOf[2]{%
    \settowidth\@TempWidthOne{$\m@th #1#2$}%
    \settowidth\@TempWidthTwo{$\m@th #1\mspace{1mu}$}%
    \strip@pt \dimexpr \@TempWidthOne*\p@/\@TempWidthTwo \relax
}
\newcommand*\SetToWidthInMathUnits[3]{%
    \settowidth\@TempWidthOne{$\m@th #2#3$}%
    \settowidth\@TempWidthTwo{$\m@th #2\mspace{1mu}$}%
    \edef #1{%
        \strip@pt \dimexpr \@TempWidthOne*\p@/\@TempWidthTwo \relax
    }%
}
\newcommand*\ExtractHexDigitOfFamilyOfSymbol[1]{%
    \expandafter\@ExtractHexDigitOfFamily \meaning #1%
}
\@ifdefinable\@ExtractHexDigitOfFamily{
    \def\@ExtractHexDigitOfFamily#1"#2#3#4#5{#3}
}
\makeatother

\newcommand*{\ReportValueOfOneEm}[1]{%
    &\text{in #1 size,}%
        &\quad\texttt{1em}%
            &=\texttt{\the \fontdimen 6 \csname #1font\endcsname 2}%
        &&\qquad\text{(from \texttt{\fontname \csname #1font\endcsname 2})}%
}
\newcommand*{\ReportSizeOfFontOfSymbol}[2]{%
    &\texttt{\pdffontsize\csname #2font\endcsname
                "\ExtractHexDigitOfFamilyOfSymbol{#1}}%
        &&\quad\text{in #2 size}%
        &\qquad(\texttt{\string\fontname}%
            &=\texttt{\fontname\csname #2font\endcsname
                "\ExtractHexDigitOfFamilyOfSymbol{#1}})%
}

\title%
[Limiting measures of supersingularities]%
{Limiting measures of supersingularities}
\author{Bodan Arsovski}
\date{December 2020}

\begin{document}
\maketitle

\newcommand{\CONS}[3]{\Lambda_{#1}({#2},{#3})}
\newcommand{\ERRRR}{\lfloor\log_p (k-1)\rfloor}
\newcommand{\ERRR}{\mathcal E}
\newcommand{\HECKEQ}[1]{T_{\mathrm{q},{#1}}}
\newcommand{\HECKES}[1]{T_{\mathrm{s},{#1}}}
\newcommand*{\myprime}{{\prime}\mkern-1.2mu}
\newcommand*{\mydprime}{{\prime\prime}\mkern-1.2mu}
\newcommand*{\mytrprime}{{\prime\prime\prime}\mkern-1.2mu}
\newcommand{\EPSISMALL}{\epsilon}
\newcommand{\ETASMALL}{\eta}
\newcommand{\DELTASMALL}{\delta}
\newcommand{\aaa}{\mathfrak{a}}
\newcommand{\COMPACTIND}[2]{\ind^{#2}}
\newcommand{\FAMILY}[1]{({#1})}
\renewcommand{\HECKET}{T}
\newcommand{\MATRIXNO}{\SCR M_1}
\newcommand{\CONSTA}{\SCR A}
\newcommand{\CONSTAMED}{10^5}
\newcommand{\CONSTAMEDMED}{10^4}
\newcommand{\CONSTB}{\SCR B}
\newcommand{\CONSTDOUBB}{2\cdot10^3}
\newcommand{\CONSTC}{10^8}
\newcommand{\ESTA}{\lfloor\log_p k\rfloor}
\newcommand{\ESTB}{\CONSTB}
\newcommand{\ESTC}{\CONSTC}
\newcommand{\MATRIXYES}{\SCR M_2}
\renewcommand{\TEMPORARYF}{m}
\renewcommand{\TEMPORARYDISC}{\SCR N}
\newcommand{\LARGESTSLOPE}{100}
\newcommand{\LARGESTPRIME}{139522597633}
\newcommand{\ddd}{\mathrm{d}}
\newcommand{\polyr}{z}
\newcommand{\polyH}{\Phi}
\newcommand{\polys}{\psi}
\newcommand{\polyconst}{C}
\newcommand{\veeC}{\check{C}}
\newcommand{\vveeC}{\hat{C}}
\newcommand{\annd}{\text{ \& }}
\newcommand{\Cnegativeone}{\frac{(-1)^{\beta}(s-\alpha)(\alpha-\beta+1)}{\beta^2(2\alpha-s+1) \Binom{\alpha}{\beta}} \Binom{\alpha}{s-\alpha}}
\newcommand{\rowR}{\mathbf{r}}
\newcommand{\colA}{\mathbf{c}}
\newcommand{\abcNL}{\\[3pt]}
\newcommand{\Cjnewline}{\\[5pt]}
\newcommand{\tempheightaer}{\rule{0pt}{12pt}}
\newcounter{equ}
\newcommand{\jj}{t}
\newcommand{\llll}{\gamma}
\newcommand{\zs}{u}
\newcommand{\za}{v}
\newcommand{\zb}{w}
\newcommand{\zi}{t}
\newcommand{\varu}{e}
\newcommand{\boxcol}[1]{\tikz{\path[draw=black,fill={#1}] (0,0) rectangle (0.2cm,0.2cm);}}
\newcommand{\RHOOO}{\varrho}\newcommand{\UUU}{U}
\newcommand{\DARKCIRC}{\bullet}
\newcommand{\subm}{sub}
\newcommand{\CONSTT}{10^6}
\newcommand{\ADJUST}{h}
\newcommand{\rr}{t}
\newcommand{\Sqbracketsss}[4]{{#1}\bullet_{%
\,#3}{#4}}
\newcommand{\Binom}[2]{\vphantom{\binom{k}{\beta}}\binom{#1}{#2}}
\newcommand{\RHOOOPRIME}{\RHOOO^\prime}

\begin{abstract}
   Let \m{p} be a prime number and let \m{k\GEQSLANT 2} be an integer. 
   In this article we study the semi-simple reductions modulo \m{p} of two-dimensional irreducible crystalline \m{p}-adic Galois representations with Hodge-Tate weights \m{0} and \m{k-1} and large slopes. Berger--Li--Zhu proved by using the theory of \m{(\varphi,\Gamma)}-modules that this reduction is constant when the slope is larger than \m{\lfloor\frac{k-2}{p-1}\rfloor}. Recently, Bergdall--Levin improved this bound to \m{\lfloor\frac{k-1}{p}\rfloor} by using the theory of Kisin modules. In this article, under the extra assumptions \m{p>3}  and \m{p+1\nmid k-1}, we asymptotically improve this bound further to \m{\lfloor\frac{k-1}{p+1}\rfloor+\ERRRR}, which is off from the predicted optimal bound \m{\approx\frac{k-1}{p+1}} only by a factor of \m{\BIGO{\log_p k}} rather than by a factor that is linear in \m{k}. As a consequence we deduce a partial result towards a conjecture by Gouv\^ea: that the measures of supersingularities of level \m{Np} oldforms tend to the zero measure on the interval \m{(\frac{1}{p+1},\frac{p}{p+1})} when \m{p} is coprime to \m{6N} and \m{\Gamma_0(N)}-regular. It is very likely that our methods extend to the cases \m{p\in\{2,3\}} and \m{p+1\nmid k-1} as well, and therefore can be adapted to eliminate the extra assumptions \m{p>3} and \m{p+1\nmid k-1}.
\end{abstract}

\section{Introduction}
\label{secInt}

\subsection{Motivation}
Let \m{p} be a prime number and let \m{k\GEQSLANT 2} and \m{N} be positive integers such that \m{N} is coprime to \m{p}. 
In~\cite{Gou01}, Gouv\^ea studied the ``supersingularities'' of weight \m{k}, level \m{\Gamma_0(Np)} eigenforms. To be specific, if \m{f} is such an eigenform with slope \m{v}, we define the supersingularity of \m{f} as \m{\frac{v}{k-1}}. The supersingularities of weight \m{k}, level \m{\Gamma_0(Np)} eigenforms belong to the interval \m{[0,1]} and are symmetric under \m{\eta \leftrightarrow 1 - \eta}, and we denote by \m{\mu_k} their discrete probability measure, \IE the probability measure on the interval \m{[0, 1]} we obtain by putting a point mass at each supersingularity.
 For \m{N=1}, his computations showed that the measure \m{\mu_k} is supported on \m{[0,\frac{1}{p+1}]\cup[\frac{p}{p+1},1]} most of the time, and the only exceptions that  occurred in his computations were for the primes \m{p} belonging to the set \Algn{\label{eqwe-0or30i349t}E_1=\{59,79,2411,3371,15271,64709,187441,27310421\}.} Even the exceptional supersingularities seemed to approach either \m{\frac{1}{p+1}} or \m{\frac{p}{p+1}} as \m{k\to\infty}, leading Gouv\^ea to the following conjecture.
 
 \begin{conjecture}[Gouv\^ea]\label{conjA}
 \it The sequence \m{(\mu_l)_{l\GEQSLANT2}} converges to the uniform measure on \m{[0,\frac{1}{p+1}]\cup[\frac{p}{p+1},1]}.
 \end{conjecture}
 
 This is the \m{p}-adic version of an interesting twist on the Sato--Tate conjecture: while the Sato--Tate conjecture asks about the distribution of the (real) slopes of a fixed modular form for varying primes, here one is interested in the distributions of the (\m{p}-adic) slopes of varying modular forms for a fixed prime.
 
 Each newform of level \m{\Gamma_0(N)} maps onto a pair of oldforms \m{f_1,f_2} of level \m{\Gamma_0(Np)} via the maps defined on \m{q}-expansions that send \m{q\mapsto q} and \m{q\mapsto q^p}. If the slope of \m{f} is \m{\alpha}, then the slopes of \m{f_1,f_2}
 are \m{\alpha} and \m{k-1-\alpha}; moreover, every newform of level \m{\Gamma_0(Np)} has slope \m{\frac{k-2}{2}} (see section~1 of~\cite{Gou01}). This is why \m{\mu_k} is symmetric under \m{\eta \leftrightarrow 1 - \eta}, and it implies that finding the slopes of eigenforms of level \m{\Gamma_0(Np)} is equivalent to  finding the slopes of newforms of level \m{\Gamma_0(N)}.

In~\cite{Buz05}, Buzzard introduced the notion of ``\m{\Gamma_0(N)}-regularity'' and made some concrete conjectures about the slopes of modular forms for \m{\Gamma_0(N)}-regular primes. This notion is related to Gouv\^ea's conjecture because all of the exceptional primes in the set \m{E_1} are \m{\Gamma_0(1)}-irregular, which led to the prediction that \m{\mu_k} is supported on \m{[0,\frac{1}{p+1}]\cup[\frac{p}{p+1},1]} when the prime \m{p} is \m{\Gamma_0(N)}-regular. One avenue to attempt to prove this prediction is via (\m{p}-adic) Galois representations, because the notion of \m{\Gamma_0(N)}-regularity can be rephrased in terms of Galois representations when \m{p>3}: a prime \m{p>3} is \m{\Gamma_0(N)}-regular if and only if, for all weights \m{l\GEQSLANT 2} and
all \m{f \in S_l (\Gamma_0(N))}, the modulo \m{p} Galois representation associated with \m{f} is reducible.
This also naturally leads to the more general prediction that all exceptions to Gouv\^ea's observation happen only for modular forms whose associated modulo \m{p} Galois representation is reducible.

 \subsection{Galois representations}
The \m{p}-adic Galois representations \m{V_{k,a}} were introduced by Colmez--Fontaine in~\cite{CF00}, and we denote their semi-simple reductions modulo \m{p} by \m{\overline V_{k,a}}. They relate to the Galois representations associated with modular forms because the Galois representation \m{\overline\rho_f} associated with an eigenform \m{f} of weight \m{k}, level \m{\Gamma_0(N)}, character \m{\chi}, and \m{U_p}-eigenvalue \m{a_p} is \Algn{\label{eqwe9i2309ri9023ri}\rho_f = {V}_{k,a_p\sqrt{\chi}}\otimes \sqrt{\chi}.} Therefore, if \m{p>3} is \m{\Gamma_0(N)}-regular then the statement 
\Algn{\label{eqe90i09ir092r}\text{``\m{\overline V_{k,a}} reducible \m{\Longrightarrow} \m{\PADICVALUATION(a)\LEQSLANT \alpha(k)} for a function \m{\alpha : \ZZ \to \QQ}''}} implies that \m{\mu_k} is supported on \m{[0,\frac{\alpha(k)}{k-1}] \cup [1-\frac{\alpha(k)}{k-1},1]}.

Berger--Li--Zhu proved in~\cite{BLZ04} by using the theory of \m{(\varphi,\Gamma)}-modules that the representations \m{\overline V_{k,a}} are locally constant around the infinite slope \m{a=0}, and they gave an explicit lower bound for the radius of the region of local constancy. To be more specific, they showed the following theorem.

\begin{theorem}[Berger--Li--Zhu]
\label{BLZ}
If \m{\PADICVALUATION(a)>\lfloor\frac{k-2}{p-1}\rfloor} then \m{\overline V_{k,a} \cong \overline V_{k,0}}.
\end{theorem}

 Since \m{\overline{V}_{k,0}\cong \ind(\omega_2^{k-1})} is irreducible when \m{k} is even (see for example proposition 6.1.2 in~\cite{Bre03b}), and there are no nontrivial eigenforms of odd weight,  theorem~\ref{BLZ} implies the following corollary.
 
\begin{corollary}\label{MAINCORO0}
If \m{p>3} is \m{\Gamma_0(N)}-regular then each of restrictions \Algn{\label{eq65f675f657f76tg678g87yh78h78}\left(\mu_l|_{(\frac1{p-1},\frac {p-2}{p-1})}\right)_{l\GEQSLANT2}} is the zero measure on \m{(\frac1{p-1},\frac {p-2}{p-1})}.
\end{corollary}

 In light of Gouv\^ea's conjecture, the lower bound \m{\lfloor\frac{k-2}{p-1}\rfloor} in theorem~\ref{BLZ} is believed to be suboptimal. The optimal bound is believed to be closer to \m{\frac{k-1}{p+1}}. Recently, Bergdall--Levin improved the bound to \m{\lfloor\frac{k-1}{p}\rfloor} by using the theory of Kisin modules (\cite{BL1}). To be more specific, they showed the following theorem.
 
 \begin{theorem}[Bergdall--Levin]
\label{BL}
If \m{\PADICVALUATION(a)>\lfloor\frac{k-1}{p}\rfloor} then \m{\overline V_{k,a} \cong \overline V_{k,0}}.
\end{theorem}
 
 A consequence of this theorem is the followning corollary, which is the analogue of corollary~\ref{MAINCORO0} with the interval \m{(\frac1{p-1},\frac {p-2}{p-1})} replaced by the interval \m{(\frac1{p},\frac {p-1}{p})}.
  
  \begin{corollary}\label{MAINCORO0prime}
If \m{p>3} is \m{\Gamma_0(N)}-regular then each of restrictions \Algn{\label{eqwe09ir09ir9304it9i39845}\left(\mu_l|_{(\frac1{p-1},\frac {p-2}{p-1})}\right)_{l\GEQSLANT2}} is the zero measure on \m{(\frac1{p},\frac {p-1}{p})}.
\end{corollary}
\subsection{Main results of this article}
For simplicity, we assume the extra condition \m{p>3}. It is very likely that our methods extend to the case \m{p\in\{2,3\}} as well, and therefore can be adapted to eliminate the extra assumption \m{p>3}.

\begin{customthm}{M}\label{MAINTHEOREM}
If \m{p>3} and \m{p+1\nmid k-1} and \Algn{\label{eq2r039ir2093tri39804tu89}\PADICVALUATION(a)>\lfloor\frac{k-1}{p+1}\rfloor+\ERRRR} then \m{\overline V_{k,a} \cong \overline V_{k,0}}.
\end{customthm}

In other words, under the extra conditions \m{p>3} and \m{p+1\nmid k-1}, we asymptotically improve the bounds in theorems~\ref{BLZ} and~\ref{BL} to \Algn{\label{eq0934it3409it4390}\lfloor\frac{k-1}{p+1}\rfloor+\ERRRR = \lfloor\frac{k-1}{p+1}\rfloor + \BIGO{\log_p k},} which is off from the the predicted optimal bound \m{\approx\frac{k-1}{p+1}} only by an additive factor of \m{\BIGO{\log_p k}} rather than by an additive factor that is linear in \m{k}. %
 
 The case \m{p+1\nmid k-1} is of local interest, but if \m{p>3} and \m{p+1\mid k-1} then \m{2\nmid k}, and there are no nontrivial eigenforms of odd weight and level \m{\Gamma_0(N)}, so theorem~\ref{MAINTHEOREM} is enough to imply the following global result.  %

\begin{customcor}{C}\label{MAINCORO}
If \m{p>3} is \m{\Gamma_0(N)}-regular then the sequence of restrictions \Algn{\label{eqwe-092ir0934it3948t}\left(\mu_l|_{(\frac1{p+1},\frac p{p+1})}\right)_{l\GEQSLANT2}} converges to the zero measure on \m{(\frac1{p+1},\frac p{p+1})}.
\end{customcor}

With conjecture~\ref{conjA} in mind, the interval \m{(\frac1{p+1},\frac p{p+1})} is optimal here.

The proof of theorem~\ref{MAINTHEOREM} is based on the local Langlands correspondence, and there are several known results that we use. We rely on Berger--Li--Zhu's theorem for the region already covered by it. Outside that region we rely on a theorem by Chenevier--Colmez on the continuity of a family of trianguline representations that \m{\overline V_{k,a}} belongs to. Counter-intuitively, the \m{\GL_2(\QQ_p)}-representation associated with \m{\overline V_{k,a}} itself via the local Langlands correspondence is too unwieldy, and the application of Chenevier--Colmez's theorem is crucial as it gives us a host of nearby representations to work with instead. We use the approximation results from~\cite{Ars1} to relate \m{\overline V_{k,a}} to these nearby representations, and thus get the best of both worlds: the simplicity of the \m{\GL_2(\QQ_p)}-representation associated with \m{\overline V_{k,a}} (a result of the ``weight'' parameter \m{k} of \m{\overline V_{k,a}} being small relative to the slope), and the structural richness of  the \m{\GL_2(\QQ_p)}-representations associated with the nearby representations (a result of the weight parameters of these nearby representations being large relative to the slope).

\subsection{Organization} %
 In section~\ref{secNot} (``\textsc{Definitions}'') we introduce all of the necessary definitions, and we also recall the \m{p}-adic and modulo \m{p} local Langlands correspondences in the form that is the most convenient for this article. In section~\ref{secLem} (``\textsc{General results about \texorpdfstring{\m{\GL_2(\QQ_p)}}{GL₂(Qₚ)}-representations}'') we prove some general results about the \m{\GL_2(\QQ_p)}-representations associated with Galois representations like \m{\overline V_{k,a}} via the local Langlands correspondence. Section~\ref{secPro} (``\textsc{Proofs}'') is the most difficult one: it contains the bulk of the ideas that make up the proof of theorem~\ref{MAINTHEOREM}. The proof contains a rather significant amount of rather difficult combinatorics. In order to improve readability, we defer many of the combinatorial results to section~\ref{secCom} (``\textsc{Combinatorics}''); the reader is encouraged to treat that section as a black box on first reading. As we use several results from~\cite{Ars1}, in order to be economical with space we refer to that article for some general and combinatorial results in sections~\ref{secLem} and~\ref{secCom}.

\subsection*{\texorpdfstring{Acknowledgements}{Acknowledgements}}\label{ackref}
I would like to thank professor Kevin Buzzard for his helpful remarks, suggestions, and support. This work was supported by Imperial College London and its President's PhD Scholarship.

\section{Definitions}
\label{secNot}

\subsection{Notation} For the remainder of this article we assume that \m{p} is a prime number and \m{k} and \m{N} are positive integers. For simplicity, we assume that \m{p>3}, though it is very likely that this assumption can be removed by suitably adjusting our combinatorial results. We also assume that \m{N\in\ZZ\backslash p\ZZ} is coprime to \m{p}. We assume that \m{a \in \ZZp} is such that \m{\PADICVALUATION(a)>\lfloor \frac{k-1}{p+1}\rfloor\GEQSLANT0}. The modulo~\m{p} Galois representation \m{\overline{V}_{k,a}} is defined in section~2 of~\cite{Ars1} and subsection~1.1 of~\cite{BLZ04}, so we do not reproduce its definition here. 
 In addition to \m{\overline{V}_{k,a}} we also consider ``nearby'' representations
\m{\overline{V}_{k^\prime,a^\prime}}, where \m{k^\prime \approx k} (\IE \m{k^\prime} is very close to \m{k} in the \m{p}-adic norm) and \m{a^\prime \approx a}. To define these precisely, we note that a local constancy theorem by Chenevier--Colmez (see proposition~4.13 in~\cite{Col08} and proposition~3.9 in~\cite{Che13} for the results, and section~4 of~\cite{Ars2} for further clarification on how they are applied in this setting) implies that \Algn{\label{eqgetryj67i8i67u567i7}\overline V_{k,a} \cong \overline V_{k+(p-1)p^m,\aaa_m}} for some integer \m{M=M(k,a)>0}, all integers \m{m \GEQSLANT M}, and a certain sequence \m{(\aaa_m)_{m\GEQSLANT M}} consisting of elements of \m{\ZZp} such that \m{\PADICVALUATION(\aaa_m)=\PADICVALUATION(a)} for all \m{m\GEQSLANT M}. For the remainder of this article we fix an integer \m{M=M(k,a)} and a sequence \m{(\aaa_m)_{m\GEQSLANT M}} with this property. We define
\Algn{\label{eqgethrytuk78o}
\nonumber r& =\textstyle k-2,\\ \textstyle \RHOOO & \textstyle = \lfloor \frac{k-1}{p+1}\rfloor = \lfloor \frac{r+1}{p+1}\rfloor,\nonumber\\ \ERRR&=\ERRRR.}
 We want to choose integers \m{\EPSISMALL,\DELTASMALL,\ETASMALL}  that  satisfy
\Algn{\label{eq498tu4958yu9458}\PADICVALUATION(\EPSISMALL) \gg \PADICVALUATION(\DELTASMALL) \gg \PADICVALUATION(\ETASMALL),} and an element \m{\aaa\in\ZZp} which has the same valuation as \m{a} such that \m{\overline V_{k+\EPSISMALL,\aaa}\cong \overline V_{k,a}}. That way we can compute the representation \m{\overline V_{k,a}} by computing the isomorphic ``nearby'' representation \m{\overline V_{k+\EPSISMALL,\aaa}}. The following explicitly chosen parameters (which depend on \m{k} and \m{a}) work:
\Algn{\label{eqr345yu9[09[90}
\nonumber\ETASMALL & \vphantom{p^{100Mk\lceil\PADICVALUATION(a)\rceil}}\textstyle = p^{100kM\lceil\PADICVALUATION(a)\rceil}\in\ZZ,\\  \nonumber\DELTASMALL &\vphantom{p^{100Mk\lceil\PADICVALUATION(a)\rceil}} \textstyle = \ETASMALL^{100\ETASMALL}\in\ZZ, \\\nonumber \EPSISMALL &\vphantom{p^{100Mk\lceil\PADICVALUATION(a)\rceil}}\textstyle = (p-1)\DELTASMALL^{100\DELTASMALL}=(p-1)p^{1000000kM\lceil\PADICVALUATION(a)\rceil \ETASMALL\DELTASMALL}\in\ZZ,\\\nonumber \rr &\vphantom{p^{100Mk\lceil\PADICVALUATION(a)\rceil}}\textstyle= r + \EPSISMALL ,\\ \aaa &\vphantom{p^{100Mk\lceil\PADICVALUATION(a)\rceil}}\textstyle = \aaa_{1000000kM\lceil\PADICVALUATION(a)\rceil \ETASMALL\DELTASMALL}\in \ZZp, \text{ satisfying } \PADICVALUATION(\aaa)=\PADICVALUATION(a) > \RHOOO=\lfloor\frac{k-1}{p+1}\rfloor,\nonumber\\ \textstyle{}\overline{V} & \textstyle{}=\vphantom{p^{100Mk\lceil\PADICVALUATION(a)\rceil}} \overline{V}_{\rr+2,\aaa}=\overline{V}_{k+\EPSISMALL,\aaa}\cong \overline{V}_{k,a} \text{ by equation~\eqref{eqgetryj67i8i67u567i7}}.}
 In light of theorem~\ref{BLZ}, theorem~\ref{MAINTHEOREM} is true if \m{\PADICVALUATION(a)>\lfloor\frac{k-2}{p-1}\rfloor}. Therefore, for the remainder of this article we assume additionally that
\Algn{\label{eq934irt934u894ut4} \lfloor\frac{k-2}{p-1}\rfloor\GEQSLANT\PADICVALUATION(a)>\lfloor\frac{k-1}{p+1}\rfloor+\ERRR.}
In particular, we note that \m{a^2\ne 4p^{k-1}}  and \m{a\ne \pm(1+p^{-1})p^{k/2}} and \m{\aaa^2\ne 4p^{k+\EPSISMALL-1}}  and \m{\aaa\ne \pm(1+p^{-1})p^{(k+\EPSISMALL)/2}}; these are eigenvalues that could potentially cause problems with the local Langlands correspondence.

\subsection{Local Langlands}
Let \m{B} be the Borel subgroup of  \m{G=\GL_2(\QQ_p)}  consisting of the upper triangular elements, 
 let
 \m{K=\GL_2(\ZZ_p) \subset G=\GL_2(\QQ_p)}, and let \m{Z} be the center of \m{G}. Let \m{\mu_x} be the unramified character of the Weil group that sends the geometric Frobenius to \m{x}, and let \m{|\ | : \QQ_p^\times \to \QQ_p^\times \hookrightarrow\QQptimes} be the \m{p}-adic norm. Let \m{W} be a finite-dimensional locally algebraic representation of the closed subgroup \m{KZ} of \m{G}. We define the compact induction of \m{W} by
    \Algn{\label{eqwr-023r-02o3r-}\textstyle \ind^G W \defeq
     \textstyle \big\{\text{locally algebraic } G \to W \,\,\big|\,\, & f (h g) = h f(g) \text{ for all } h \in KZ \nonumber \\ &\textstyle{} \text{\& }\Supp f\text{ is compact in }KZ\backslash G \}.}
 Suppose that \m{W} is over the field \m{\FF \in \{\QQp,\FFp\}}. For 
  elements \m{g\in G} and \m{w \in W} we write \m{\sqbrackets{g}{H}{\FF}{w}} for the unique element of \m{\ind^GW} that is supported on \m{KZg^{-1}} and maps \m{g^{-1}} to \m{w}. Every element of \m{\ind^G W} can be written as a finite linear combination of functions of the type \m{\sqbrackets{g}{H}{\FF}{w}}, and \Algn{\label{eqwrei2093ri029ir902ri}\textstyle g_1(\sqbrackets{g_2}{H}{\FF}{(h w)}) = \sqbrackets{(g_1g_2h)}{H}{\FF}{w}.}
 For \m{l\GEQSLANT0} we define
\Algn{\label{eqr098wu4t9845ut894}\textstyle\TILDE\SIGMA_{l} = \textstyle{\uline{\Sym}^{l}}({\overline{\QQ}_p^2)} \defeq {\Sym^{l}}({\overline{\QQ}_p^2)} \otimes \absdet^{l/2},} and we define \m{\SIGMA_{l}} as the reduction  of \m{\uline{\Sym}^{l}({\overline{\ZZ}_p^2)}} modulo  the  maximal ideal \m{\maxideal} of \m{\ZZp}. For \m{h\in\ZZ} we define \Algn{\label{eqwe90i2309i209r}\sigma_h \defeq{} \uline{\Sym}^h(\FFp^2).}
As in section~3 of~\cite{Ars1}, we note that we can view these as \m{G}-modules of homogeneous polynomials in two variables.

Let us also define the Hecke operator \m{T \in \End_G(\COMPACTIND{KZ}{G} \TILDE\SIGMA_{\rr})} corresponding to the double coset of \m{\begin{Smatrix}
p & 0 \\ 0 & 1
\end{Smatrix}}. This operator satisfies the explicit formula  
    \Algn{ \label{eq3r984ti394tu98t}& \textstyle \HECKET ({\gamma}\Sqbracketsss{}{KZ}{\QQp}{v}) %
     \textstyle = {\sum_{\MU\in\FF_p} \gamma\begin{Lmatrix}
    p &[\MU] \\ 0 & 1
    \end{Lmatrix} }\Sqbracketsss{}{KZ}{\QQp}{\left(\begin{Lmatrix}
    1 & -[\MU] \\ 0 & p
    \end{Lmatrix} \cdot v\right)} + \gamma\begin{Lmatrix}
    1 & 0 \\ 0 & p
    \end{Lmatrix} \Sqbracketsss{}{KZ}{\QQp}{\left(\begin{Lmatrix}
    p & 0 \\ 0 & 1
    \end{Lmatrix} \cdot v\right)},} where \m{[\TEMPORARYA]} is the Teichm\"uller lift of \m{\TEMPORARYA \in \FF_p} to \m{\ZZ_p}. 
 We define 
\Algn{\label{eqf9384jut984ut9845}\nonumber\textstyle \Pi_{\rr+2,\aaa} & \textstyle= \COMPACTIND{KZ}{G} \TILDE\SIGMA_{\rr}/(T-a),\\\textstyle \Theta_{\rr+2,\aaa} & \nonumber\textstyle = \Img\left(\COMPACTIND{KZ}{G}(\uline{\Sym}^{\rr}(\overline{\ZZ}_p^2)) \xrightarrow{\hspace{\SIZEEE} \hphantom{\cong} \hspace{\SIZEEE}}\Pi_{\rr+2,\aaa}\right), \\
\textstyle{}\overline\Theta_{\rr+2,\aaa} & \textstyle = \Theta_{\rr+2,\aaa} \otimes_{\ZZp} \FFp.} In particular, \m{\overline{\Theta}_{\rr+2,\aaa}} is a quotient of \m{\COMPACTIND{KZ}{G} \SIGMA_{\rr}}, and we define \m{\SCR I} to be the ideal such that \Algn{\label{eqfo43093i09r4}\overline{\Theta}_{\rr+2,\aaa} \cong \COMPACTIND{KZ}{G} \SIGMA_{\rr}/\SCR I.}
 The ideal \m{\SCR I} contains the reduction modulo~\m{p} of any integral element in the image of \m{T-\aaa}. For \m{j\in\{0,\ldots,p-1\}}, \m{\lambda \in\FFp}, and a character \m{\psi:\QQ_p^\times \to \overline\FF_p^\times}, we write %
    \Algn{\label{eq09i439tui3498t5}\textstyle \repnn(j,\lambda,\psi) \defeq{} ({\COMPACTIND{KZ}{G} \sigma_t}/{(T_{\sigma}-\lambda)})\otimes\psi,} where \m{T_\sigma\in \End_G(\COMPACTIND{KZ}{G} \sigma_j)} is the Hecke operator corresponding to the double coset of \m{\begin{Smatrix} p&0 \\ 0 &1\end{Smatrix}}. We let \m{\omega} be the modulo~\m{p} reduction of the cyclotomic character, \m{\ind(\omega_2^{j+1})} be the unique irreducible representation whose determinant is \m{{\omega^{j+1}}} and that is equal to \m{{\omega_2^{j+1}\oplus\omega_2^{p(j+1)}}} on inertia, \m{\MORE{h}\in\{1,\ldots,p-1\}} and \m{\LESS{h}\in\{0,\ldots,p-2\}} be the numbers in the corresponding sets that are congruent to \m{h} modulo \m{p-1}.
The following theorem is the main result of \cite{Ber10} and says that the modulo \m{p} Langlands correspondence is compatible with the \m{p}-adic local Langlands correspondence.
 
\begin{theorem}\label{BERGERX1}
          There are \m{j\in\{0,\ldots,p-1\}} and \m{\psi:\QQ_p^\times \to \overline\FF_p^\times} such that either 
          \Algn{\label{eqwe09or3094it9034ti}\textstyle \overline\Theta_{k,a} \cong \repnn(j,0,\psi) %
          }
          or
          \Algn{\label{eqwr2or09r3094t}\textstyle \overline\Theta_{k,a}^{\rm ss} \cong \left(\repnn(j,\lambda,\psi) \oplus\repnn(\LESS{p-3-j},\lambda^{-1},\omega^{j+1}\psi) \right)^{\rm ss}}
          for some \m{\lambda \in\FFp}.
          In the former case we have
          \Algn{\label{eqw3rpoir09ri3490}\textstyle \overline V_{k,a} \cong \ind(\omega_2^{j+1})\otimes\psi,}
          and in the latter case we have
          \Algn{\label{eqwe092ir09ri3940}\textstyle \overline V_{k,a} \cong \left(\mu_\lambda \omega^{j+1} \oplus \mu_{\lambda^{-1}}\right)\otimes\psi.}
   \end{theorem}

Let \m{\overline\Theta = \overline \Theta_{\rr+2,\aaa}^{\rm ss}}. Proposition~4.1.4 in~\cite{BLZ04} implies that \Algn{\label{eqi4309i30} \overline V_{k,0} \cong \left\{\begin{array}{rl}
     \ind(\omega_2^{k-1}) & \text{if } p+1 \nmid k-1\text{,}   \\[5pt]
     \left(\mu_{\sqrt{-1}}\oplus \mu_{-\sqrt{-1}}\right)\otimes \omega^{(k-1)/(p+1)} & \text{if } p+1\mid k-1\text{.} 
 \end{array}\right.} Therefore, in order to prove theorem~\ref{MAINTHEOREM}, we want to show that 
 \Algn{\label{eqdrftgyuhihygt}\overline V = \overline V_{\rr+2,\aaa} \cong  \ind(\omega_2^{k-1}),
}
(since \m{p+1\nmid k-1}). 
So theorem~\ref{BERGERX1} implies that theorem~\ref{MAINTHEOREM} can be rewritten in the following equivalent form.

\begin{customthm}{M'}\label{MAINTHEOREMp}
Recall that \m{p>3} and \m{p+1\nmid k-1} and \m{\PADICVALUATION(\aaa)>\lfloor\frac{k-1}{p+1}\rfloor+\ERRR}. We have
\Algn{\label{qowur983u4t3yt8374t4y}\overline \Theta  \cong
     \repnn(\overline{r-2\RHOOO},0,\omega^\RHOOO) .}
\end{customthm}

So the goal of this article is to prove theorem~\ref{MAINTHEOREMp}.

\subsection{More notation} Let \m{\TEMPORARYI=\lfloor\PADICVALUATION(\aaa)\rfloor+1}. We cite  section~4 of~\cite{Ars1} for the definitions of \m{\BIGO{\alpha}}, and (for \m{h\in\ZZ}) \m{I_h}, and \m{\theta=xy^p-x^py}, and the evaluation \m{[P]} of a boolean \m{P}, which we do not reproduce here. We also recall that there is a filtration
\Algn{\label{eqwr0929u398ut} \overline\Theta = \overline\Theta_0 \supset \overline\Theta_1 \supset \cdots \supset \overline \Theta_\alpha \supset \cdots \supset \overline{\Theta}_\TEMPORARYI=0}
whose \m{\alpha}th subquotient (for \m{\alpha\in\{0,\ldots,\TEMPORARYI-1\}}) is a subquotient of \Algn{\label{eqe09i2riw90ri} \FILT N_\alpha=\ind^G\left(\overline{\theta}^{\alpha}\Sigma_{\rr-\alpha(p+1)}/\overline{\theta}^{\alpha+1}\Sigma_{\rr-(\alpha+1)(p+1)}\right) \cong \ind^G I_{\rr-2\alpha}(\alpha) \cong \ind^G I_{r-2\alpha}(\alpha).}
 To be specific, if an element of \Algn{\label{eqwer092ir948}\ind^G\left(\overline{\theta}^{\alpha}\Sigma_{\rr-\alpha(p+1)}/\overline{\theta}^{\alpha+1}\Sigma_{\rr-(\alpha+1)(p+1)}\right)} is represented by an element of \m{\SCR I \subset \COMPACTIND{KZ}{G}\SIGMA_\rr}, then that element is trivial in the subquotient \m{\overline \Theta_\alpha/\overline\Theta_{\alpha+1}}. Finally, we define (for \m{\alpha\in\{0,\ldots,\TEMPORARYI-1\}})
\Algn{\label{eqw90eri2309ri3904}\nonumber
    \SUBMODULE(\alpha) & \textstyle = \sigma_{\MORE{r-2\alpha}}(\alpha) \subset N_\alpha,\\
    \QUOTIENTI(\alpha) & = N_\alpha/\sigma_{\MORE{r-2\alpha}}(\alpha) \cong \sigma_{\LESS{2\alpha-r}}(r-\alpha),
} similarly as in section~4 of~\cite{Ars1}, and we denote by \m{\HECKEQ{\alpha},\HECKES{\alpha}} the Hecke operators corresponding to the double coset of \m{\begin{Smatrix} p&0 \\ 0 &1\end{Smatrix}} on the modules \m{ \COMPACTIND{KZ}{G}\QUOTIENTI(\alpha),\COMPACTIND{KZ}{G}\SUBMODULE(\alpha)}, respectively.
For \m{\alpha\in\{0,\ldots,\DELTASMALL\}} we define
\Algn{\label{eqwrwt34t34y45y453}\ADJUST_\alpha & \textstyle=  x^{\alpha}y^{\rr-\alpha}-x^{\alpha+\DELTASMALL}y^{\rr-\alpha-\DELTASMALL}\in\TILDE\SIGMA_\rr,\nonumber\\ \textstyle\ADJUST_\alpha^\ast & =\textstyle \begin{Lmatrix}
      0 & 1\\ 1 & 0      \end{Lmatrix} \ADJUST_\alpha=x^{\rr-\alpha}y^{\alpha}-x^{\rr-\alpha-\DELTASMALL}y^{\alpha+\DELTASMALL}\in\TILDE\SIGMA_\rr.}
 For \m{\alpha,\beta,R\GEQSLANT0} we define \m{\CONS{R}{\alpha}{\beta}} by
\Algn{\label{eq0943if0943if0934i} \sum_{\beta=\alpha-R}^\alpha \CONS{R}{\alpha}{\beta} \Binom{(p-1)X+\alpha}{\alpha-\beta} = \Binom{R-X}{R} = (-1)^R \Binom{X-1}{R} \in \QQ_p[X].}
Note that both sides of equation~\eqref{eq0943if0943if0934i} are polynomials in \m{X} over \m{\QQ_p} of degree \m{R}.

\section{General results about \texorpdfstring{\m{\GL_2(\QQ_p)}}{GL₂(Qₚ)}-representations}
\label{secLem}

\begin{lemma}\label{LEMMAX2}
If \m{\alpha \in \{0,\ldots,\ETASMALL\}} then
  \Algn{\label{eqf093409if039i30}
      \textstyle & \textstyle \aaa \Sqbracketsss{}{KZ}{\QQp}{ \ADJUST_\alpha} 
      \EQUIVZEROIDEAL  p^{\alpha} \begin{Lmatrix}
      1 & 0 \\ 0 & p
      \end{Lmatrix} \Sqbracketsss{}{KZ}{\QQp}{x^{\alpha}y^{\rr-\alpha}} 
   + \BIGO{p^{2\ETASMALL}}.
  }
If \m{\alpha \in \{0,\ldots,\ETASMALL\}}, \m{\beta\in\{\alpha,\ldots,\ETASMALL\}}, and \m{\FAMILY{C_l}_{l\in \ZZ}} is a family of elements of \m{\ZZ_p} then
\Algn{\label{eq0r3i4i30309ri4309ir}\nonumber
      & \textstyle \sum_{i} \left(\sum_{l=\alpha-\beta}^{\alpha} C_l \Binom{\rr-\alpha+l}{i(p-1)+l}\right)\Sqbracketsss{}{KZ}{\QQp}{x^{i(p-1)+\alpha}y^{\rr-i(p-1)-\alpha}} \\ & \textstyle\null\qquad \EQUIVZEROIDEAL \frac{\aaa p^{-\alpha}}{p-1} \sum_{l=\alpha-\beta}^{\alpha} C_{l}p^l\sum_{\MU\in\FF_p^\times} [\MU]^{-l}\begin{Lmatrix}
p & [\MU] \\ 0 & 1
\end{Lmatrix}  \Sqbracketsss{}{KZ}{\QQp}{\ADJUST_{\alpha-l}}+\BIGO{p^{\ETASMALL}}.}
\end{lemma}

\begin{proof}{Proof. }
 We have
\Algn{\label{eqf938498ugf3984}
\aaa  \Sqbracketsss{}{KZ}{\QQp}{  \ADJUST_\alpha} & \textstyle \EQUIVZEROIDEAL T\left(1  \Sqbracketsss{}{KZ}{\QQp}{ \ADJUST_\alpha} \right) \nonumber\\ &\textstyle 
\EQUIVZEROIDEAL {\sum_{\MU\in\FF_p} \begin{Lmatrix}
    p &[\MU] \\ 0 & 1
    \end{Lmatrix} }\Sqbracketsss{}{KZ}{\QQp}{A_\MU} + \begin{Lmatrix}
    1 & 0 \\ 0 & p
    \end{Lmatrix} \Sqbracketsss{}{KZ}{\QQp}{A},
}
where, due to the explicit equation for \m{T} (equation~\eqref{eq3r984ti394tu98t}), 
\renewcommand{\tempheightaer}{\vphantom{(-[\MU])^{\MORE{\rr-\alpha-\TEMPORARYA}} \sum_{j,\TEMPORARYA\GEQSLANT 0}   (-1)^{j} \Binom{\ETASMALL}{j}\Binom{\rr-j(p-1)-\alpha}{\TEMPORARYA}p^\TEMPORARYA x^{\rr-\TEMPORARYA }y^\TEMPORARYA}}\Algn{\label{eqf3uj498ut4398ut489u} A_\MU & \textstyle \tempheightaer=\nonumber 
x^{\alpha}(-[\MU]x+py)^{\rr-\alpha}-x^{\alpha+\DELTASMALL}(-[\MU]x+py)^{\rr-\alpha-\DELTASMALL}
 \\ & \textstyle
\tempheightaer=  \sum_{\TEMPORARYA\GEQSLANT 0} (-[\MU])^{\MORE{\rr-\alpha-\TEMPORARYA}} \left(\Binom{\rr-\alpha}{\TEMPORARYA} -\Binom{\rr-\alpha-\DELTASMALL}{\TEMPORARYA} \right) p^\TEMPORARYA x^{\rr-\TEMPORARYA }y^\TEMPORARYA = \BIGO{\DELTASMALL p^{-2\ETASMALL} + p^{2\ETASMALL}} = \BIGO{p^{2\ETASMALL}},
} 
and \Algn{\label{eqr2098u4f398ug4398ujr9348} A & \textstyle = p^\alpha x^\alpha y^{\rr-\alpha} + \BIGO{p^{\alpha+\DELTASMALL}} = p^\alpha x^\alpha y^{\rr-\alpha} + \BIGO{p^{2\ETASMALL}}.
}

Equations~\eqref{eqf938498ugf3984}, \eqref{eqf3uj498ut4398ut489u}, and~\eqref{eqr2098u4f398ug4398ujr9348} imply equation~\eqref{eqf093409if039i30}. Equation~\eqref{eqf093409if039i30} implies that
\Algn{\label{eqe09rtgtr23ir0943ir430r}\nonumber
& \textstyle \frac{\aaa p^{-\alpha}}{p-1} \sum_{l=\alpha-\beta}^\alpha C_{l}p^l\sum_{\MU\in\FF_p^\times} [\MU]^{-l}\begin{Lmatrix}
p & [\MU] \\ 0 & 1
\end{Lmatrix}  \Sqbracketsss{}{KZ}{\QQp}{ \ADJUST_{\alpha-l}}\nonumber 
\\ \nonumber&\textstyle \null\qquad \EQUIVZEROIDEAL \frac{1}{p-1} \sum_{l=\alpha-\beta}^\alpha C_{l}\sum_{\MU\in\FF_p^\times} [\MU]^{-l}\begin{Lmatrix}
1 & [\MU] \\ 0 & 1
\end{Lmatrix}  \Sqbracketsss{}{KZ}{\QQp}{ x^{\alpha-l}y^{\rr-\alpha+l}} + \BIGO{p^{2\ETASMALL-\beta}}
\\ &\textstyle \null\qquad \EQUIVZEROIDEAL  \sum_{l=\alpha-\beta}^\alpha C_{l} \sum_i \Binom{\rr-\alpha+l}{i(p-1)+l}
\Sqbracketsss{}{KZ}{\QQp}{ x^{i(p-1)+\alpha}y^{\rr-i(p-1)-\alpha}}
+ \BIGO{p^{\ETASMALL}},} which implies equation~\eqref{eq0r3i4i30309ri4309ir}.
\end{proof}

\begin{lemma}\label{LEMMAX3}
{
Let \m{\alpha \in \{0,\ldots,\ETASMALL\}} and \m{v \in \QQ} and the family \m{\FAMILY{D_i}_{i\in \ZZ}} of elements of \m{\ZZ_p} be such that
\Algn{\label{eqfew94if09i094if4309f}\nonumber\textstyle %
D_i&\textstyle{}=0 \text{ for } i\not\in[\frac{-\alpha}{p-1},\frac{\rr-\alpha}{p-1}],%
\\\nonumber v 
& \textstyle  \LEQSLANT \PADICVALUATION(\TEMPORARYD_{w}(D_\DARKCIRC)) \text{ for } \alpha \LEQSLANT w \LEQSLANT 2\ETASMALL,
\\ v %
 & \textstyle <\PADICVALUATION(\TEMPORARYD_{w}(D_\DARKCIRC)) \text{ for } 0\LEQSLANT w<\alpha.
} %
For \m{j \in \ZZ}, let
\Algn{\label{eq0493i093irf3}\textstyle\DELTAA_j = (-1)^{j-\ETASMALL} (1-p)^{-\alpha} \Binom{\alpha}{j-\ETASMALL} \TEMPORARYD_{\alpha}(D_\DARKCIRC),}
so that \m{\FAMILY{\DELTAA_j}_{j\in \ZZ}} is supported on the set of indices \m{\{\ETASMALL,\ldots,\alpha+\ETASMALL\}} and therefore \m{\TEMPORARYD_{w}(\DELTAA_\DARKCIRC)} is properly defined for \m{0\LEQSLANT w < \alpha}.
Then \m{v \LEQSLANT\PADICVALUATION(\TEMPORARYD_{\alpha}(\DELTAA_\DARKCIRC))\LEQSLANT \PADICVALUATION(\DELTAA_j)} for all \m{j \in \ZZ}, and
\renewcommand{\tempheightaer}{\vphantom{ \sum_{i \in \left[0,(2\PADICVALUATION(\aaa)-\alpha)/(p-1)\right] \cap \ZZ}  \left(D_{i}
 \Sqbracketsss{}{KZ}{\QQp}{\ADJUST_{i(p-1)+\alpha}}+ D_{\rr-i}
 \Sqbracketsss{}{KZ}{\QQp}{\ADJUST^\ast_{i(p-1)+\alpha}}\right)}}
\Algn{\label{equ3984u4398ur3498ur}\nonumber
      \textstyle & \textstyle
       \sum_i (\DELTAA_i - D_i) \Sqbracketsss{}{KZ}{\QQp}{x^{i(p-1)+\alpha}y^{\rr-i(p-1)-\alpha}}
      \\ & \textstyle \null\qquad 
      \tempheightaer\EQUIVZEROIDEAL
       - \sum_{i \LEQSLANT (\ETASMALL+\PADICVALUATION(\aaa)-\alpha)/(p-1)}  D_{i}
 \Sqbracketsss{}{KZ}{\QQp}{\ADJUST_{i(p-1)+\alpha}}\nonumber
 \\ & \textstyle \null\phantom{\EQUIVZEROIDEAL}\qquad\qquad   \tempheightaer
       - \sum_{ i\GEQSLANT (\rr-\alpha-\ETASMALL-\PADICVALUATION(\aaa))/(p-1)}  D_{i}
 \Sqbracketsss{}{KZ}{\QQp}{ 
\ADJUST^\ast_{\rr-i(p-1)-\alpha} }\nonumber
  \\ & \textstyle \null\phantom{\EQUIVZEROIDEAL}\qquad\qquad   \tempheightaer+ E \Sqbracketsss{}{KZ}{\QQp}{\theta^{\alpha+1}h} + %
  F \Sqbracketsss{}{KZ}{\QQp}{h^{\prime}} + \BIGO{p^\ETASMALL},
  }
  for some polynomials %
 \m{h,h^{\prime}} and some %
  \m{E,F\in\ZZp} with %
  \m{\PADICVALUATION(E)\GEQSLANT v} and \m{\PADICVALUATION(F)>v}.
}
\end{lemma}

\begin{proof}{Proof. }
By using the equation
\Algn{\label{eqr09i03i4t093i4t} {\gamma}\Sqbracketsss{}{KZ}{\QQp}{v} \EQUIVZEROIDEAL \aaa^{-1} \HECKET ({\gamma}\Sqbracketsss{}{KZ}{\QQp}{v}), }
and equation~\eqref{eq3r984ti394tu98t}, we can deduce that
\Algn{\label{eq34j983ut9834u}\nonumber& \textstyle\tempheightaer
\sum_i (\DELTAA_i - D_i) \Sqbracketsss{}{KZ}{\QQp}{x^{i(p-1)+\alpha}y^{\rr-i(p-1)-\alpha}}
\\ 
& \textstyle\null\qquad\tempheightaer \EQUIVZEROIDEAL \aaa^{-1}T \left(\sum_i (\DELTAA_i - D_i) \Sqbracketsss{}{KZ}{\QQp}{x^{i(p-1)+\alpha}y^{\rr-i(p-1)-\alpha}}\right)\nonumber
\\ 
& \textstyle\null\qquad \EQUIVZEROIDEAL\tempheightaer \aaa^{-1}\sum_i (\DELTAA_i - D_i) \sum_{\LAMBDA\in\FF_p^\times } \begin{Lmatrix} p & [\LAMBDA] \\ 0 & 1
 \end{Lmatrix}\Sqbracketsss{}{KZ}{\QQp}{x^{i(p-1)+\alpha}(-[\LAMBDA]x+py)^{\rr-i(p-1)-\alpha}} \nonumber
\\
& \nonumber\textstyle\null\phantom{\EQUIVZEROIDEAL}\tempheightaer\qquad \qquad + \aaa^{-1}\sum_{i}(\DELTAA_i - D_i)  \left( p^{\rr-i(p-1)-\alpha}\begin{Lmatrix} p & 0 \\ 0 & 1
 \end{Lmatrix} +p^{i(p-1)+\alpha}\begin{Lmatrix} 1 & 0 \\ 0 & p
 \end{Lmatrix}\right)
 \nonumber\\ & \textstyle\null\qquad\qquad\qquad\qquad\qquad\qquad\qquad\qquad\qquad\qquad\qquad\tempheightaer\Sqbracketsss{}{KZ}{\QQp}{x^{i(p-1)+\alpha}y^{\rr-i(p-1)-\alpha}}\nonumber
\\ 
& \textstyle\null\qquad\tempheightaer \EQUIVZEROIDEAL \aaa^{-1}\sum_i (\DELTAA_i - D_i) \sum_{\LAMBDA\in\FF_p^\times } \begin{Lmatrix} p & [\LAMBDA] \\ 0 & 1
 \end{Lmatrix}\Sqbracketsss{}{KZ}{\QQp}{x^{i(p-1)+\alpha}(-[\LAMBDA]x+py)^{\rr-i(p-1)-\alpha}} \nonumber
   \\ & \textstyle \null\phantom{\EQUIVZEROIDEAL}\qquad \qquad
      \tempheightaer 
       - \sum_{i \LEQSLANT (\ETASMALL+\PADICVALUATION(\aaa)-\alpha)/(p-1)}  D_{i}
 \Sqbracketsss{}{KZ}{\QQp}{\ADJUST_{i(p-1)+\alpha}}\nonumber
 \\ & \textstyle \null\phantom{\EQUIVZEROIDEAL}\qquad\qquad   \tempheightaer
       - \sum_{ i\GEQSLANT (\rr-\alpha-\ETASMALL-\PADICVALUATION(\aaa))/(p-1)}  D_{i}
 \Sqbracketsss{}{KZ}{\QQp}{ 
\ADJUST^\ast_{\rr-i(p-1)-\alpha} } + \BIGO{p^\ETASMALL}.
} The third congruence follows from lemma~\ref{LEMMAX2}. We also have
\Algn{\label{eq958gu4u55}\nonumber & \textstyle{}\tempheightaer \sum_{\LAMBDA\in\FF_p^\times } \begin{Lmatrix} p & [\LAMBDA] \\ 0 & 1
 \end{Lmatrix}\Sqbracketsss{}{KZ}{\QQp}{x^{i(p-1)+\alpha}(-[\LAMBDA]x+py)^{\rr-i(p-1)-\alpha}}\\ &\null\tempheightaer\qquad\textstyle \EQUIVZEROIDEAL \sum_{\TEMPORARYA=0}^{2\ETASMALL}  \Binom{\rr-i(p-1)-\alpha}{\TEMPORARYA} p^{\TEMPORARYA} \sum_{\LAMBDA\in\FF_p^\times } [-\LAMBDA]^{\rr-\alpha-\TEMPORARYA}\begin{Lmatrix} p & [\LAMBDA] \\ 0 & 1
 \end{Lmatrix} \Sqbracketsss{}{KZ}{\QQp}{x^{\rr-\TEMPORARYA} y^{\TEMPORARYA}}+\BIGO{p^{2\ETASMALL}}\nonumber
 \\ &\null\tempheightaer\qquad\textstyle \EQUIVZEROIDEAL \aaa \sum_{\TEMPORARYA=0}^{2\ETASMALL}  \Binom{\rr-i(p-1)-\alpha}{\TEMPORARYA}   \Sqbracketsss{}{KZ}{\QQp}{\sum_{\LAMBDA\in\FF_p^\times } [-\LAMBDA]^{\rr-\alpha-\TEMPORARYA}\begin{Lmatrix} 1 & [\LAMBDA] \\ 0 & 1
 \end{Lmatrix}\ADJUST^\ast_{\TEMPORARYA}}+\BIGO{p^{2\ETASMALL}}.}
 The second congruence follows from lemma~\ref{LEMMAX2}. By assumption, if \Algn{\label{eqf90i4f3409if30if309i}X_\TEMPORARYA =  \sum_i (\DELTAA_i - D_i) \Binom{\rr-i(p-1)-\alpha}{\TEMPORARYA},} then \m{\PADICVALUATION(X_\TEMPORARYA)>v} for \m{\TEMPORARYA\in\{0,\ldots,\alpha\}}, and \m{\PADICVALUATION(X_\TEMPORARYA)\GEQSLANT v}  for \m{\TEMPORARYA\in\{\alpha+1,\ldots,2\ETASMALL\}}. This means that equation~\eqref{eq958gu4u55} implies that 
 \Algn{\label{eq4093i309i309ir4}\nonumber&\textstyle \aaa^{-1}\sum_i (\DELTAA_i - D_i) \sum_{\LAMBDA\in\FF_p^\times } \begin{Lmatrix} p & [\LAMBDA] \\ 0 & 1
 \end{Lmatrix}\Sqbracketsss{}{KZ}{\QQp}{x^{i(p-1)+\alpha}(-[\LAMBDA]x+py)^{\rr-i(p-1)-\alpha}} \\ & \null\tempheightaer\qquad\textstyle\EQUIVZEROIDEAL \sum_{\TEMPORARYA=0}^{2\ETASMALL}  X_\TEMPORARYA \Sqbracketsss{}{KZ}{\QQp}{\sum_{\LAMBDA\in\FF_p^\times } [-\LAMBDA]^{\rr-\alpha-\TEMPORARYA}\begin{Lmatrix} 1 & [\LAMBDA] \\ 0 & 1
 \end{Lmatrix} \ADJUST^\ast_{\TEMPORARYA}}+\BIGO{p^{\ETASMALL}}, }
 which together with equation~\eqref{eq34j983ut9834u} implies equation~\eqref{equ3984u4398ur3498ur} with
 \Algn{\label{eq43i0ri309ir3409}\nonumber E\theta^{\alpha+1}h &\textstyle\tempheightaer= \sum_{\TEMPORARYA=\alpha+1}^{2\ETASMALL}  X_\TEMPORARYA  \sum_{\LAMBDA\in\FF_p^\times } [-\LAMBDA]^{\rr-\alpha-\TEMPORARYA}\begin{Lmatrix} 1 & [\LAMBDA] \\ 0 & 1
 \end{Lmatrix} \ADJUST^\ast_{\TEMPORARYA},\\\textstyle Fh^\prime &\textstyle\tempheightaer= \sum_{\TEMPORARYA=0}^{\alpha}  X_\TEMPORARYA  \sum_{\LAMBDA\in\FF_p^\times } [-\LAMBDA]^{\rr-\alpha-\TEMPORARYA}\begin{Lmatrix} 1 & [\LAMBDA] \\ 0 & 1
 \end{Lmatrix} \ADJUST^\ast_{\TEMPORARYA}.}%
\end{proof}

\begin{lemma}\label{LEMMAX4}
{
Let  \m{\FAMILY{C_l}_{l\in \ZZ}} be any family of elements of \m{\ZZ_p}. Suppose that 
\m{\alpha \in \{0,\ldots, \ETASMALL\}} and \m{\beta\in\{\alpha,\ldots,\ETASMALL\}} and \m{v \in \QQ} and   the family \m{\FAMILY{D_i}_{i\in \ZZ}} defined by
 \Algn{\label{eqf4309i309i4309}\textstyle D_i = %
 \VERIFYIF{i\in\{\lceil\frac{-\alpha}{p-1}\rceil,\ldots,\lfloor\frac{\rr-\alpha}{p-1}\rfloor\}}D_i^\prime + \sum_{l=\alpha-\beta}^\alpha C_{l} \Binom{\rr-\alpha+l}{i(p-1)+l}}
satisfy %
\Algn{\label{eq09i3209id3094if4309fi}\textstyle  \textstyle v  & \textstyle  \LEQSLANT \PADICVALUATION(\TEMPORARYD_{w}(D_\DARKCIRC)) \text{ for } \alpha \LEQSLANT w \LEQSLANT 2\DELTASMALL,
\nonumber\\ v & \textstyle <\PADICVALUATION(\TEMPORARYD_{w}(D_\DARKCIRC)) \text{ for } 0\LEQSLANT w<\alpha.
} Note that \m{\FAMILY{D_i}_{i\in \ZZ}} is supported on the finite set of indices \m{\{\lceil\frac{-\alpha}{p-1}\rceil,\ldots,\lfloor\frac{\rr-\alpha}{p-1}\rfloor\}}. %
 Then
\Algn{\label{eqf943i903i3r}\nonumber
      \textstyle & \textstyle
       (1-p)^{-\alpha} \TEMPORARYD_{\alpha}(D_\DARKCIRC)  \Sqbracketsss{}{KZ}{\QQp}{\theta^\alpha x^{\ETASMALL(p-1)}y^{\rr-\alpha(p+1)-\ETASMALL(p-1)}}
       \\ & \textstyle \null\qquad \nonumber
      \EQUIVZEROIDEAL \frac{\aaa p^{-\alpha}}{p-1} \sum_{l=\alpha-\beta}^{\alpha} C_{l}p^l\sum_{\MU\in\FF_p^\times} [\MU]^{-l}\begin{Lmatrix}
p & [\MU] \\ 0 & 1
\end{Lmatrix}  \Sqbracketsss{}{KZ}{\QQp}{\ADJUST_{\alpha-l}}
\\ & \textstyle \null\phantom{\EQUIVZEROIDEAL}\qquad\qquad 
      \tempheightaer
+ \sum_{i} D_i^\prime\Sqbracketsss{}{KZ}{\QQp}{x^{i(p-1)+\alpha}y^{\rr-i(p-1)-\alpha}}\nonumber
      \\ & \textstyle \null\phantom{\EQUIVZEROIDEAL}\qquad\qquad 
      \tempheightaer
       - \sum_{i \LEQSLANT (\ETASMALL+\PADICVALUATION(\aaa)-\alpha)/(p-1)}  D_{i}
 \Sqbracketsss{}{KZ}{\QQp}{\ADJUST_{i(p-1)+\alpha}}\nonumber
 \\ & \textstyle \null\phantom{\EQUIVZEROIDEAL}\qquad\qquad   \tempheightaer
       - \sum_{ i\GEQSLANT (\rr-\alpha-\ETASMALL-\PADICVALUATION(\aaa))/(p-1)}  D_{i}
 \Sqbracketsss{}{KZ}{\QQp}{ 
 \ADJUST^\ast_{\rr-i(p-1)-\alpha} }\nonumber
  \\ & \textstyle \null\phantom{\EQUIVZEROIDEAL}\qquad\qquad   \tempheightaer+ E \Sqbracketsss{}{KZ}{\QQp}{\theta^{\alpha+1}h} + %
  F \Sqbracketsss{}{KZ}{\QQp}{h^{\prime}} + \BIGO{p^\ETASMALL},
  }
  for some polynomials %
 \m{h,h^{\prime}} and some %
  \m{E,F\in\ZZp} with %
  \m{\PADICVALUATION(E)\GEQSLANT v} and \m{\PADICVALUATION(F)>v}.

}
\end{lemma}

\begin{proof}
{Proof. }
Lemma~\ref{LEMMAX2} implies that
\Algn{\label{eqf439fi349i98}\nonumber
      & \textstyle \sum_{i} (D_i-D_i^\prime)\Sqbracketsss{}{KZ}{\QQp}{x^{i(p-1)+\alpha}y^{\rr-i(p-1)-\alpha}} \\ & \textstyle\null\qquad \EQUIVZEROIDEAL \frac{\aaa p^{-\alpha}}{p-1} \sum_{l=\alpha-\beta}^{\alpha} C_{l}p^l\sum_{\MU\in\FF_p^\times} [\MU]^{-l}\begin{Lmatrix}
p & [\MU] \\ 0 & 1
\end{Lmatrix}  \Sqbracketsss{}{KZ}{\QQp}{\ADJUST_{\alpha-l}}+\BIGO{p^{\ETASMALL}}.} Equation~\eqref{eqf439fi349i98} together with  lemma~\ref{LEMMAX3} implies that
\Algn{\label{eqf43if948u5f98}\nonumber
      \textstyle & \textstyle
       \sum_i \DELTAA_i \Sqbracketsss{}{KZ}{\QQp}{x^{i(p-1)+\alpha}y^{\rr-i(p-1)-\alpha}}
       \\ & \textstyle \null\qquad \nonumber
      \EQUIVZEROIDEAL \frac{\aaa p^{-\alpha}}{p-1} \sum_{l=\alpha-\beta}^{\alpha} C_{l}p^l\sum_{\MU\in\FF_p^\times} [\MU]^{-l}\begin{Lmatrix}
p & [\MU] \\ 0 & 1
\end{Lmatrix}  \Sqbracketsss{}{KZ}{\QQp}{\ADJUST_{\alpha-l}} \\ & \textstyle \null\phantom{\EQUIVZEROIDEAL}\qquad\qquad 
      \tempheightaer
+ \sum_{i} D_i^\prime\Sqbracketsss{}{KZ}{\QQp}{x^{i(p-1)+\alpha}y^{\rr-i(p-1)-\alpha}}\nonumber
      \\ & \textstyle \null\phantom{\EQUIVZEROIDEAL}\qquad\qquad 
      \tempheightaer
       - \sum_{i \LEQSLANT (\ETASMALL+\PADICVALUATION(\aaa)-\alpha)/(p-1)}  D_{i}
 \Sqbracketsss{}{KZ}{\QQp}{\ADJUST_{i(p-1)+\alpha}}\nonumber
 \\ & \textstyle \null\phantom{\EQUIVZEROIDEAL}\qquad\qquad   \tempheightaer
       - \sum_{ i\GEQSLANT (\rr-\alpha-\ETASMALL-\PADICVALUATION(\aaa))/(p-1)}  D_{i}
 \Sqbracketsss{}{KZ}{\QQp}{ 
 \ADJUST^\ast_{\rr-i(p-1)-\alpha} }\nonumber
  \\ & \textstyle \null\phantom{\EQUIVZEROIDEAL}\qquad\qquad   \tempheightaer+ E \Sqbracketsss{}{KZ}{\QQp}{\theta^{\alpha+1}h} + %
  F \Sqbracketsss{}{KZ}{\QQp}{h^{\prime}} + \BIGO{p^\ETASMALL},
  }
  for some polynomials %
 \m{h,h^{\prime}} and some %
  \m{E,F\in\ZZp} with %
  \m{\PADICVALUATION(E)\GEQSLANT v} and \m{\PADICVALUATION(F)>v}. Equation~\eqref{eqf43if948u5f98} can be rewritten in the form of equation~\eqref{eqf943i903i3r} because
  \Algn{\label{eqf-04i3r9i390i}\nonumber & \textstyle\sum_i \DELTAA_i \Sqbracketsss{}{KZ}{\QQp}{x^{i(p-1)+\alpha}y^{\rr-i(p-1)-\alpha}} \\ & \textstyle \null\qquad = (1-p)^{-\alpha} \TEMPORARYD_{\alpha}(D_\DARKCIRC)  \Sqbracketsss{}{KZ}{\QQp}{\theta^\alpha x^{\ETASMALL(p-1)}y^{\rr-\alpha(p+1)-\ETASMALL(p-1)}}.}%
\end{proof}

\section{Proofs}
\label{secPro}

 We want to prove theorem~\ref{MAINTHEOREMp} by computing \m{\overline\Theta}. We accomplish this as the cumulative result of the following six subsections.

\subsection{\protect\boldmath If \texorpdfstring{\m{Q}}{Q} is an \texorpdfstring{\m{\infty}}{∞}-dimensional factor of \texorpdfstring{\m{\overline{\Theta}}}{\textTheta}, then \texorpdfstring{\m{Q}}{Q} is not a factor of \texorpdfstring{\m{\FILT N_\alpha}}{N\textalpha}, for \texorpdfstring{\m{\alpha\in\{0,\ldots,\RHOOO-1\}}}{\textalpha\ in 0,...,\textrho-1}}
 \label{subsecProf1}
For
\m{\alpha\in\{0,\ldots,\RHOOO-1\}}, let us define the matrix
\Algn{\label{eq09i435t98ut498r9384ut3948} M_\alpha^{(r)} = \left(\Binom{r-\alpha+j}{i(p-1)+j}\right)_{\{i \,|\,  i(p-1)+\alpha \in(\RHOOO,r-\RHOOO)\},\,\alpha-\RHOOO\LEQSLANT j\LEQSLANT\alpha}.} So the rows of 
 \m{M_\alpha^{(r)} } are indexed by those \m{i} such that \m{i(p-1)+\alpha \in(\RHOOO,r-\RHOOO)}, and the columns of \m{M_\alpha^{(r)} } are indexed by those \m{j} such that \m{\alpha-\RHOOO\LEQSLANT j\LEQSLANT\alpha}. This means that \m{M_\alpha^{(r)} } has \m{C=\RHOOO+1} columns  and
\Algn{\label{eqf304i093i309i} R \LEQSLANT \lfloor\frac{r-2\rho+p-2}{p-1}\rfloor \LEQSLANT \RHOOO+1= C
}
rows, \IE \m{M_\alpha^{(r)} } has no more rows than columns. Let
\Algn{\label{eq093ri290ir}
M_\alpha^{(r)\prime}=\left(\Binom{r-\alpha+j}{i(p-1)+j}\right)_{\{i \,|\,  i(p-1)+\alpha \in(\RHOOO,r-\RHOOO)\},\,\alpha-R< j\LEQSLANT\alpha}}
be the right \m{R\times R} submatrix of \m{M_\alpha^{(r)}}. We can write
\Algn{\label{eq093i40fi4}\Binom{r-\alpha+j}{i(p-1)+j} = \Binom{r}{i(p-1)+\alpha}\Binom{i(p-1)+\alpha}{\alpha-j}\Binom{r}{\alpha-j}^{-1},} so the \m{\ZZ_p}-module determined by the image of the matrix \m{M_\alpha^{(r)\prime}} contains the \m{\ZZ_p}-module determined by the image of the matrix \Algn{\label{eqwr9i34t09i9u4589uy9845u9y}\left(\Binom{r}{i(p-1)+\alpha}\Binom{i(p-1)+\alpha}{\alpha-j}\right)_{\{i \,|\,  i(p-1)+\alpha \in(\RHOOO,r-\RHOOO)\},\,\alpha-R< j\LEQSLANT\alpha}.}
Lemma~\ref{lemmaD} implies that \Algn{\label{eqwproi230r9i3049ti3} \PADICVALUATION\left(\Binom{r}{i(p-1)+\alpha}\right) \LEQSLANT \lfloor\log_p(r+1)\rfloor=\ERRR,} so the latter \m{\ZZ_p}-module  contains \m{p^{\ERRR}} \m{\times} the 
\m{\ZZ_p}-module determined by the image of the matrix
\Algn{\label{eq409ir093ir3209r}M_\alpha^{(r)\mydprime} = \left(\Binom{ i(p-1)+\alpha}{\alpha-j}\right)_{\{i \,|\,  i(p-1)+\alpha \in(\RHOOO,r-\RHOOO)\},\,\alpha-R< j\LEQSLANT\alpha}.}
There exists a \m{\gamma\in\ZZ_{\GEQSLANT0}} such that \m{M_\alpha^{(r)\mydprime}} is obtained from
\Algn{\label{eq90ir0394ri4039}{M_\alpha^{(r)\mytrprime}} = \left(\Binom{i(p-1)+\gamma}{j}\right)_{0\LEQSLANT i,j<R}} by permuting the rows.
By Vandermonde's convolution formula,
\Algn{\label{eq049309r349i309}M_\alpha^{(r)\mytrprime} = \left(\Binom{i(p-1)}{j}\right)_{0\LEQSLANT i,j<R} \cdot \left(\Binom{\gamma}{j-i}\right)_{0\LEQSLANT i,j<R}.}
 Since the matrix \Algn{\label{eqwr-02it3409t498ut439}\left(\Binom{\gamma}{j-i}\right)_{0\LEQSLANT i,j<R}} is upper triangular with \m{1}'s on the diagonal and \Algn{\label{eq903i430iir4309ri30ti30}\det\left(\Binom{i(p-1)}{j}\right)_{0\LEQSLANT i,j<R} = (p-1)^R \det\left(\Binom{i}{j}\right)_{0\LEQSLANT i,j<R}= (p-1)^R} by a variant of Vandermonde's determinant identity, the reduction modulo~\m{p} of \m{M_\alpha^{(r)\mytrprime}} has full rank (in characteristic~\m{p}). This in turn implies that the reduction modulo~\m{p} of \m{M_\alpha^{(r)\mydprime}} has full rank.
 Therefore, for each \m{u} such that \Algn{\label{eqwe-09oi2309ri34it934t} u(p-1) +\alpha\in (\RHOOO,r-\RHOOO),} there exist constants \m{C_\alpha(r,u),\ldots,C_{\alpha-R+1}(r,u)} such that
 \Algn{\label{eq094ti309tiu9485t}M_\alpha^{(r)\prime} \left(C_\alpha(r,u),\ldots,C_{\alpha-R+1}(r,u)\right)^T = p^{\ERRR}(\VERIFYIF{i=u})_{\{i \,|\,  i(p-1)+\alpha \in(\RHOOO,r-\RHOOO)\}},}
 \IE such that
\Algn{\label{eqf0934i0439i4309}\sum_{l=\alpha-R+1}^\alpha C_l(r,u) \Binom{r-\alpha+l}{i(p-1)+l}= \VERIFYIF{i=u} p^{\ERRR}} for all \m{i} such that \Algn{\label{eqe092ir3498rt349tu4598t}i(p-1)+\alpha \in(\RHOOO,r-\RHOOO).} By adding linear combinations of equation~\eqref{eqf0934i0439i4309} for varying \m{u}, we get that 
\Algn{\label{eqf34ik093i0}\nonumber &\textstyle \sum_{ i(p-1)+\alpha \in(\RHOOO,r-\RHOOO)} \sum_{l=\alpha-\RHOOO}^\alpha C_l \Binom{r-\alpha+l}{i(p-1)+l} x^{i(p-1)+\alpha}y^{r-i(p-1)-\alpha}  \\ & \null\qquad\textstyle +\sum_{i(p-1)+\alpha\in [0,\RHOOO]\cup[r-\RHOOO,r]} D_i^\prime x^{i(p-1)+\alpha}y^{r-i(p-1)-\alpha}\nonumber \\& \null\qquad\qquad\textstyle = p^{\ERRR}\theta^\alpha x^{p-1}y^{r-\alpha(p+1)-p+1}\in\TILDE\SIGMA_r,} for some \m{C_l,D_i^\prime}.
Let \m{D_i(r)} be the coefficient of \m{x^{i(p-1)+\alpha}y^{r-i(p-1)-\alpha}} on the right side of equation~\eqref{eqf34ik093i0}. Then, due to part~(5) of lemma~6 and lemma~7 in~\cite{Ars1},
\Algn{\label{eq0934i093ir4309}\TEMPORARYD_w(D_\bullet(r)) =\sum_i D_i(r) \Binom{i(p-1)}{w}}
is zero for \m{0\LEQSLANT w <\alpha}, and has valuation that is greater than or equal to \m{ \ERRR} for \m{w\GEQSLANT \alpha}, with equality for \m{w=\alpha}. Let 
\m{D_i} be the coefficient of \m{x^{i(p-1)+\alpha}y^{r-i(p-1)-\alpha}}
in
\Algn{\label{eq039fi409fi409gi54}\nonumber &\textstyle \sum_{i(p-1)+\alpha \in(\RHOOO,\rr-\RHOOO)} \sum_{l=\alpha-\RHOOO}^\alpha C_l \Binom{\rr-\alpha+l}{i(p-1)+l} x^{i(p-1)+\alpha}y^{\rr-i(p-1)-\alpha}  \\ & \null\qquad\textstyle +\sum_{i(p-1)+\alpha\in [0,\RHOOO]\cup[\rr-\RHOOO,\rr]} D_i^\mydprime x^{i(p-1)+\alpha}y^{\rr-i(p-1)-\alpha},}
where \m{D_i^\mydprime = D_i^\prime} and \m{D_{\rr-i}^\mydprime = D_{r-i}^\prime} for \m{i(p-1)+\alpha\in [0,\RHOOO]}. Since \Algn{\label{eq0943i093093i0349i}\PADICVALUATION((i(p-1)+\alpha)!) \LEQSLANT \PADICVALUATION(\RHOOO!) \LEQSLANT k} for \m{i(p-1)+\alpha\in [0,\RHOOO]}, it is easy to show by using lemma~5 in~\cite{Ars1} that
\Algn{\label{eqf394i094i4095} \TEMPORARYD_w(D_\bullet) = \TEMPORARYD_w(D_\bullet(r)) + {\BIGO{\EPSISMALL p^{-k-W}}}} for all \m{0\LEQSLANT w \LEQSLANT W}. In particular, we can apply lemma~\ref{LEMMAX4} to the constants \m{(D_i)_{i\in\ZZ}} and to \m{v=\ERRR}, and as a result get that
\Algn{\label{eqf34if409i4gi409i}\nonumber
      \textstyle & \textstyle
       (1-p)^{-\alpha} \TEMPORARYD_{\alpha}(D_\DARKCIRC)  \Sqbracketsss{}{KZ}{\QQp}{\theta^\alpha x^{\ETASMALL(p-1)}y^{\rr-\alpha(p+1)-\ETASMALL(p-1)}}
       \\ & \textstyle \null\qquad \nonumber
      \EQUIVZEROIDEAL \frac{\aaa p^{-\alpha}}{p-1} \sum_{l=\alpha-\RHOOO}^{\alpha} C_{l}p^l\sum_{\MU\in\FF_p^\times} [\MU]^{-l}\begin{Lmatrix}
p & [\MU] \\ 0 & 1
\end{Lmatrix}  \Sqbracketsss{}{KZ}{\QQp}{\ADJUST_{\alpha-l}}
\\ & \textstyle \null\phantom{\EQUIVZEROIDEAL}\qquad\qquad 
      \tempheightaer
+ \sum_{i(p-1)+\alpha\in [0,\RHOOO]\cup[\rr-\RHOOO,\rr]} D_i^\mytrprime\Sqbracketsss{}{KZ}{\QQp}{x^{i(p-1)+\alpha}y^{\rr-i(p-1)-\alpha}}\nonumber
      \\ & \textstyle \null\phantom{\EQUIVZEROIDEAL}\qquad\qquad 
      \tempheightaer
       - \sum_{i \LEQSLANT (\ETASMALL+\PADICVALUATION(\aaa)-\alpha)/(p-1)}  D_{i}
 \Sqbracketsss{}{KZ}{\QQp}{\ADJUST_{i(p-1)+\alpha}}\nonumber
 \\ & \textstyle \null\phantom{\EQUIVZEROIDEAL}\qquad\qquad   \tempheightaer
       - \sum_{ i\GEQSLANT (\rr-\alpha-\ETASMALL-\PADICVALUATION(\aaa))/(p-1)}  D_{i}
 \Sqbracketsss{}{KZ}{\QQp}{ 
 \ADJUST^\ast_{\rr-i(p-1)-\alpha} }\nonumber
  \\ & \textstyle \null\phantom{\EQUIVZEROIDEAL}\qquad\qquad   \tempheightaer+ E \Sqbracketsss{}{KZ}{\QQp}{\theta^{\alpha+1}h} + %
  F \Sqbracketsss{}{KZ}{\QQp}{h^{\prime}} + \BIGO{p^\ETASMALL},
  }
  for some %
 \m{h,h^{\prime}} and some %
  \m{E,F\in\ZZp} with %
  \m{\PADICVALUATION(E)\GEQSLANT \ERRR} and \m{\PADICVALUATION(F)>\ERRR}. Here
  \Algn{\label{eq9034i03ir94ut598u4}D_i^\mytrprime = D_i^\mydprime - \sum_{l=\alpha-\RHOOO}^\alpha C_l \Binom{\rr-\alpha+l}{i(p-1)+l} } for all \m{i} such that \Algn{\label{eqwe092ri3490ti439ti} i(p-1)+\alpha\in [0,\RHOOO]\cup[\rr-\RHOOO,\rr].} The left side of equation~\eqref{eqf34if409i4gi409i} is \m{p^{\ERRR}\psi}, where \m{\psi} is an integral element whose reduction modulo~\m{p} represents a generator of \m{\FILT N_\alpha}. We can use lemma~\ref{LEMMAX2} to get that the first and second lines on the right side of equation~\eqref{eqf34if409i4gi409i} are \Algn{\label{eq98u349f8u4398uf}\BIGO{p^{\PADICVALUATION(\aaa)-\RHOOO}}=\BIGO{p^{\ERRR}}.}
 We can also use equation~\eqref{eqf34ik093i0} to get that \Algn{\label{eq90i0943i4309}D_i = \BIGO{p^{\ERRR}}} for all \m{i} such that 
 \Algn{\label{eqwe09i2r0t9i34t9i3049ti}i(p-1)+\alpha \in [0,r-\RHOOO) \cup (\rr-r+\RHOOO,\rr].} This is because
 \Algn{\label{eqwf094i0439ir3}D_i=D_i(r)+\BIGO{\EPSISMALL p^{-\PADICVALUATION(\ETASMALL!)}}} for all \m{i} such that \m{i(p-1)+\alpha \in [0,r-\RHOOO)}, and \Algn{\label{eq34ikg50953i3095}D_{\rr-i} = D_{r-i}(r) + \BIGO{\EPSISMALL p^{-\PADICVALUATION(\ETASMALL!)}}} for all \m{i} such that
 \m{i(p-1)+\alpha \in (\rr-r+\RHOOO,\rr]}.
 We also have, due to equation~\eqref{eqf34ik093i0},  \Algn{\label{eqfo3k4039ir0439ir0}D_w = \BIGO{\EPSISMALL p^{-\PADICVALUATION(\ETASMALL!)}}} for all \m{i} and all \Algn{\label{eqwe-02or093ir9034it390t5445}w\in \ZZ_{\LEQSLANT 0}\cup\ZZ_{\GEQSLANT \frac{\rr-2\alpha}{p-1}}.} This, together with lemma~\ref{LEMMAX2}, implies that the sum of the third and fourth lines on the right side of equation~\eqref{eqf34if409i4gi409i} is \m{p^{\ERRR}} \m{\times} an integral element whose reduction modulo~\m{p} represents the trivial element of \m{\FILT N_\alpha}. Finally, the fifth line of on the right side of equation~\eqref{eqf34if409i4gi409i} is evidently \m{p^{\ERRR}} \m{\times} an integral element whose reduction modulo~\m{p} represents the trivial element of \m{\FILT N_\alpha}. 
 Therefore equation~\eqref{eqf34if409i4gi409i} gives an element in \m{\SCR I} that generates  \m{\FILT N_\alpha}, implying that no \texorpdfstring{\m{\infty}}{∞}-dimensional factor of \m{\overline{\Theta}}  is a subquotient of {\m{\FILT N_\alpha}}.

 \subsection{\protect\boldmath If \texorpdfstring{\m{Q}}{Q} is an \texorpdfstring{\m{\infty}}{∞}-dimensional factor of \texorpdfstring{\m{\overline{\Theta}}}{\textTheta}, then \texorpdfstring{\m{Q}}{Q} is not a factor of \texorpdfstring{\m{\FILT N_\alpha}}{N\textalpha} for \texorpdfstring{\m{\alpha\in\{\RHOOO+1,\ldots,\TEMPORARYI-1\}}}{\textalpha\ in \textrho+1,...,\textnu-1}}
 \label{subsecProf2}
 Let
 \Algn{\label{eq2309ri3290tiu3490}\textstyle\RHOOOPRIME=\lceil\frac{r-\alpha}{p}\rceil-1.}
 For \m{\alpha\in\{\RHOOO+1,\ldots,\TEMPORARYI-1\}} and \m{l \in \{\alpha-\RHOOOPRIME,\ldots,\alpha\}} let us define \Algn{\label{eq4309it0934i}C_l = \CONS{\RHOOOPRIME}{\alpha}{l}\Binom{r}{\alpha-l}.}
 As in subsection~\ref{subsecProf1} we can conclude that
   \Algn{\label{eq9034i093ti0439it09} &\textstyle \sum_{ i(p-1)+\alpha \in (\alpha,\RHOOOPRIME(p-1)+\alpha]} \sum_{l=\alpha-\RHOOOPRIME}^\alpha C_l \Binom{r-\alpha+l}{i(p-1)+l} x^{i(p-1)+\alpha}y^{r-i(p-1)-\alpha}  = 0\in\TILDE\SIGMA_r.} %
     Let 
\m{D_i} be the coefficient of \m{x^{i(p-1)+\alpha}y^{\rr-i(p-1)-\alpha}}
in
\Algn{\label{eqf0934i09i309f430}  &\textstyle \sum_{ i(p-1)+\alpha \in (\alpha,\rr-r+\RHOOOPRIME(p-1)+\alpha]} \sum_{l=\alpha-\RHOOOPRIME}^{ \alpha} C_l \Binom{\rr-\alpha+l}{i(p-1)+l} x^{i(p-1)+\alpha}y^{\rr-i(p-1)-\alpha}. %
}
 Then it is easy to show by using lemma~5 in~\cite{Ars1} that
\Algn{\label{eq90i409t5i4t0i454} \TEMPORARYD_w(D_\bullet) = {\BIGO{\EPSISMALL p^{-k-W}}}} for all \m{0\LEQSLANT w \LEQSLANT W}. In particular, we can apply lemma~\ref{LEMMAX4} to the constants \m{(D_i)_{i\in\ZZ}} and to \m{v=\ETASMALL}, and as a result get that
\Algn{\label{eqf904095gi45i0g5i049g}\nonumber
      \textstyle & \textstyle
        \nonumber
      \frac{\aaa p^{-\alpha}}{p-1} \sum_{l=\alpha-\RHOOOPRIME}^{\alpha} C_{l}p^l\sum_{\MU\in\FF_p^\times} [\MU]^{-l}\begin{Lmatrix}
p & [\MU] \\ 0 & 1
\end{Lmatrix}  \Sqbracketsss{}{KZ}{\QQp}{\ADJUST_{\alpha-l}}
\\ & \textstyle \null\qquad {\EQUIVZEROIDEAL}
      \tempheightaer
\sum_{i(p-1)+\alpha\in[0,\alpha]\cup(\rr-r+\RHOOOPRIME(p-1)+\alpha,\rr]} D_i^\prime\Sqbracketsss{}{KZ}{\QQp}{x^{i(p-1)+\alpha}y^{\rr-i(p-1)-\alpha}}\nonumber
      \\ & \textstyle \null\phantom{\EQUIVZEROIDEAL}\qquad\qquad 
      \tempheightaer
       + \sum_{i \LEQSLANT (\ETASMALL+\PADICVALUATION(\aaa)-\alpha)/(p-1)}  D_{i}
 \Sqbracketsss{}{KZ}{\QQp}{\ADJUST_{i(p-1)+\alpha}}\nonumber
 \\ & \textstyle \null\phantom{\EQUIVZEROIDEAL}\qquad\qquad   \tempheightaer
       + \sum_{ i\GEQSLANT (\rr-\alpha-\ETASMALL-\PADICVALUATION(\aaa))/(p-1)}  D_{i}
 \Sqbracketsss{}{KZ}{\QQp}{ 
 \ADJUST^\ast_{\rr-i(p-1)-\alpha} }
    + \BIGO{p^\ETASMALL},
  }
  where
  \Algn{\label{eq094i039i}D_i^\prime =  \sum_{l=\alpha-\RHOOOPRIME}^{\alpha} C_l \Binom{\rr-\alpha+l}{i(p-1)+l} } for all \m{i} such that \Algn{\label{eqwe-0ro34t0i4398ti4958ti4985yui45} i(p-1)+\alpha\in[0,\alpha]\cup(\rr-r+\RHOOOPRIME(p-1)+\alpha,\rr].} 
By approximating \m{D_i} with \m{D_i(r)+\BIGO{\EPSISMALL p^{-\PADICVALUATION(\ETASMALL!)}}} as in subsection~\ref{subsecProf1} we can show that the third and fourth lines of equation~\eqref{eqf904095gi45i0g5i049g} are in \Algn{\label{eq0943i309i039r}\BIGO{p^{\alpha p-\PADICVALUATION(\aaa)}}\nonumber&=\textstyle\BIGO{\aaa p^{-\alpha+(k-2\PADICVALUATION(\aaa)-p+3)}}\\ & \nonumber\textstyle=\BIGO{\aaa p^{-\alpha+(p-3)(k/(p-1)-1)}}\\ & \textstyle=\BIGO{\aaa p^{-\alpha+2\ERRR}}.}
Consequently we get that
\Algn{\label{eqiu45jiu54jiu4tji4u}\nonumber
      \textstyle & \textstyle
        \nonumber
      \sum_{i(p-1)+\alpha\in[0,\alpha]\cup(\rr-r+\RHOOOPRIME(p-1)+\alpha,\rr]} \sum_{l=\alpha-\RHOOOPRIME}^{\alpha} C_l \Binom{\rr-\alpha+l}{i(p-1)+l} \Sqbracketsss{}{KZ}{\QQp}{x^{i(p-1)+\alpha}y^{\rr-i(p-1)-\alpha}}
\\ & \textstyle \null\qquad {\EQUIVZEROIDEAL}
     \frac{\aaa p^{-\alpha}}{p-1} \sum_{l=\alpha-\RHOOOPRIME}^{\alpha} C_{l}p^l\sum_{\MU\in\FF_p^\times} [\MU]^{-l}\begin{Lmatrix}
p & [\MU] \\ 0 & 1
\end{Lmatrix}  \Sqbracketsss{}{KZ}{\QQp}{\ADJUST_{\alpha-l}} \tempheightaer
   + \BIGO{\aaa p^{-\alpha+2\ERRR}}.
  }
  Lemma~\ref{LEMMAX2} and  the definition of \m{(C_l)_{\alpha-\RHOOOPRIME\LEQSLANT l\LEQSLANT\alpha}} then imply that
  \Algn{\label{eq45tik49ti09i409it45}\nonumber
      \textstyle & \textstyle
        \nonumber
      \sum_{i(p-1)+\alpha\in[0,\alpha]} X_i \begin{Lmatrix}
p & 0 \\ 0 & 1
\end{Lmatrix}  \Sqbracketsss{}{KZ}{\QQp}{\ADJUST_{i(p-1)+\alpha}}
\\ & \textstyle\null\qquad +  \sum_{i(p-1)+\alpha\in(\rr-r+\RHOOOPRIME(p-1)+\alpha,\rr]} X_i^\ast\begin{Lmatrix}
1 & 0 \\ 0 & p
\end{Lmatrix} \Sqbracketsss{}{KZ}{\QQp}{\ADJUST^\ast_{\rr-i(p-1)-\alpha}}
\nonumber\\ & \textstyle \null\qquad\qquad\phantom{+} {\EQUIVZEROIDEAL}
     \frac{1}{p-1} \sum_{l=\alpha-\RHOOOPRIME}^{\alpha} C_{l}p^l\sum_{\MU\in\FF_p^\times} [\MU]^{-l}\begin{Lmatrix}
p & [\MU] \\ 0 & 1
\end{Lmatrix}  \Sqbracketsss{}{KZ}{\QQp}{\ADJUST_{\alpha-l}} \tempheightaer
   + \BIGO{p^{2\ERRR}},
  }
  where
  \Algn{\label{eqr3094u94t5y78}\nonumber X_i & \textstyle = p^{-i(p-1)}\Binom{\rr}{i(p-1)+\alpha} \Binom{\RHOOOPRIME-i}{\RHOOOPRIME},\\ X_i^\ast & \textstyle = p^{i(p-1)+2\alpha-\rr}\Binom{\rr}{i(p-1)+\alpha} \Binom{\RHOOOPRIME-i}{\RHOOOPRIME} .}
  By lemma~\ref{lemmaC1a},
  \Algn{\label{eqofw8ueow3}\PADICVALUATION(X_0) = \PADICVALUATION\left(\Binom{t}{\alpha}\right) < \PADICVALUATION\left(p^{-i(p-1)}\Binom{\rr}{i(p-1)+\alpha} \Binom{\RHOOOPRIME-i}{\RHOOOPRIME}\right)=\PADICVALUATION(X_i)}
  for all \m{i} such that \m{i(p-1)+\alpha\in[0,\alpha)}. By lemma~\ref{lemmaC1b}, 
  \Algn{\label{eqf3j49ujt847gt87} \PADICVALUATION(X_0) =\PADICVALUATION\left(\Binom{t}{\alpha}\right) <\PADICVALUATION\left(p^{i(p-1)+2\alpha-\rr}\Binom{\rr}{i(p-1)+\alpha} \Binom{\RHOOOPRIME-i}{\RHOOOPRIME}\right) = \PADICVALUATION(X_i^\ast)}
  for all \m{i} such that \m{i(p-1)+\alpha\in(\rr-r+\RHOOOPRIME(p-1)+\alpha,\rr]}. By lemma~\ref{lemmaC2a}, 
  \Algn{\label{eqf98wuerfeferg} \PADICVALUATION(X_0) =\PADICVALUATION\left(\Binom{t}{\alpha}\right) <\PADICVALUATION\left(\CONS{\RHOOOPRIME}{\alpha}{l}\Binom{\rr}{\alpha-l}p^l\right)=\PADICVALUATION\left(C_lp^l\right)}
  for all \m{l \in \{\alpha-\RHOOOPRIME,\ldots,\alpha\}}. Moreover, by lemma~\ref{lemmaD},
  \Algn{\label{eqjf3948g849gu498} \PADICVALUATION(X_0) =\PADICVALUATION\left(\Binom{t}{\alpha}\right)=\PADICVALUATION\left(\Binom{r}{\alpha}\right) < 2\ERRR.}
  So if we divide both sides of equation~\eqref{eq45tik49ti09i409it45} by \m{\Binom{t}{\alpha}} we get an integral element, and if we reduce that integral element modulo~\m{p} then the only contributing term to the result is the ``\m{i=0}'' term in the first line of equation~\eqref{eq45tik49ti09i409it45}. %
   Therefore we can conclude that \m{\SCR I} contains
\Algn{\label{equ54t98u54e49rgtje} \begin{Lmatrix}
p & 0 \\ 0 & 1
\end{Lmatrix}  \Sqbracketsss{}{KZ}{\FFp}{\ADJUST_{\alpha}},} 
 which  
  represents a generator of \m{\FILT N_\alpha}, and we can conclude the desired result.

 \subsection{\protect\boldmath If \texorpdfstring{\m{r-\RHOOO(p+1)=p-2}}{r-\textrho(p+1)=p-2} and \texorpdfstring{\m{Q}}{Q} is an \texorpdfstring{\m{\infty}}{∞}-dimensional factor of \texorpdfstring{\m{\overline{\Theta}}}{\textTheta}, then \texorpdfstring{\m{Q}}{Q} is not a factor of \texorpdfstring{\m{\COMPACTIND{KZ}{G}\SUBMODULE(\RHOOO)}}{ind(sub(\textrho))}}
 \label{subsecProf3}

 Let \m{r-\RHOOO(p+1)=p-2}, so that \Algn{\label{eq9043r309ir}\SUBMODULE(\RHOOO)\cong \sigma_{p-2}(\RHOOO).}  The proof in this case is very similar to the proof in subsection~\ref{subsecProf1}, so we just give a rough sketch.
 We let
 \Algn{\label{eq34ik039r4309i} M^{(r)} = \left(\Binom{r-\RHOOO+j}{i(p-1)+j}\right)_{\{i \,|\,  i(p-1)+\RHOOO \in(\RHOOO,r-\RHOOO)\},\,0\LEQSLANT j\LEQSLANT\RHOOO}.} 
 As in  subsection~\ref{subsecProf1} we can prove that the image of a certain lattice under the right square submatrix of \m{M^{(r)} } (seen as an endomorphism) contains \m{p^{\ERRR}} \m{\times} that lattice.
 We can conclude the following analogous equation to  equation~\eqref{eqf34ik093i0}:
 \Algn{\label{eq984ut9u544ut}\nonumber &\textstyle \sum_{i(p-1)+\RHOOO \in(\RHOOO,r-\RHOOO)} \sum_{l=0}^\RHOOO C_l \Binom{r-\RHOOO+l}{i(p-1)+l} x^{i(p-1)+\RHOOO}y^{r-i(p-1)-\RHOOO}  \\ & \null\qquad\textstyle +\sum_{i(p-1)+\RHOOO\in [0,\RHOOO]\cup[r-\RHOOO]} D_i^\prime x^{i(p-1)+\RHOOO}y^{r-i(p-1)-\RHOOO}\nonumber \\& \null\qquad\qquad\textstyle = p^{\ERRR}\theta^\RHOOO  y^{r-\RHOOO(p+1)}\in\TILDE\SIGMA_r,} for some integers \m{D_i^\prime}.
 The main difference with equation~\eqref{eqf34ik093i0} is that we must write \m{\theta^\RHOOO y^{r-\RHOOO(p+1)}} instead of \m{\theta^\RHOOO x^{p-1}y^{r-\RHOOO(p+1)-p+1}}.
 This means that we can only conclude that \Algn{\label{eq34r903ir039i}D_w = \BIGO{\EPSISMALL p^{-\PADICVALUATION(\ETASMALL!)}}}  for all \Algn{\label{eqwer-0o4t9i4598t9485t945}w\in \ZZ_{< 0}\cup\ZZ_{\GEQSLANT \frac{\rr-2\RHOOO}{p-1}}} (rather than for all \m{w\in \ZZ_{\LEQSLANT 0}\cup\ZZ_{\GEQSLANT \frac{\rr-2\RHOOO}{p-1}}}) in the equation
 \Algn{\label{eqj45gj498ug98u}\nonumber
      \textstyle & \textstyle
       (1-p)^{-\RHOOO} \TEMPORARYD_{\RHOOO}(D_\DARKCIRC)  \Sqbracketsss{}{KZ}{\QQp}{\theta^\RHOOO x^{\ETASMALL(p-1)}y^{\rr-\RHOOO(p+1)-\ETASMALL(p-1)}}
       \\ & \textstyle \null\qquad \nonumber
      \EQUIVZEROIDEAL \frac{\aaa p^{-\RHOOO}}{p-1} \sum_{l=0}^{\RHOOO} C_{l}p^l\sum_{\MU\in\FF_p^\times} [\MU]^{-l}\begin{Lmatrix}
p & [\MU] \\ 0 & 1
\end{Lmatrix}  \Sqbracketsss{}{KZ}{\QQp}{\ADJUST_{\RHOOO-l}}
\\ & \textstyle \null\phantom{\EQUIVZEROIDEAL}\qquad\qquad 
      \tempheightaer
+ \sum_{i(p-1)+\RHOOO\in [0,\RHOOO]\cup[\rr-\RHOOO,\rr]} D_i^\mytrprime\Sqbracketsss{}{KZ}{\QQp}{x^{i(p-1)+\RHOOO}y^{\rr-i(p-1)-\RHOOO}}\nonumber
      \\ & \textstyle \null\phantom{\EQUIVZEROIDEAL}\qquad\qquad 
      \tempheightaer
       - \sum_{i \LEQSLANT (\ETASMALL+\PADICVALUATION(\aaa)-\RHOOO)/(p-1)}  D_{i}
 \Sqbracketsss{}{KZ}{\QQp}{\ADJUST_{i(p-1)+\RHOOO}}\nonumber
 \\ & \textstyle \null\phantom{\EQUIVZEROIDEAL}\qquad\qquad   \tempheightaer
       - \sum_{ i\GEQSLANT (\rr-\RHOOO-\ETASMALL-\PADICVALUATION(\aaa))/(p-1)}  D_{i}
 \Sqbracketsss{}{KZ}{\QQp}{ 
\ADJUST^\ast_{\rr-i(p-1)-\RHOOO} }\nonumber
  \\ & \textstyle \null\phantom{\EQUIVZEROIDEAL}\qquad\qquad   \tempheightaer+ E \Sqbracketsss{}{KZ}{\QQp}{\theta^{\RHOOO+1}h} + %
  F \Sqbracketsss{}{KZ}{\QQp}{h^{\prime}} + \BIGO{p^\ETASMALL},
  }
  for some %
 \m{h,h^{\prime}} and some %
  \m{E,F\in\ZZp} with %
  \m{\PADICVALUATION(E)\GEQSLANT \ERRR} and \m{\PADICVALUATION(F)>\ERRR}, which is the analogous equation to equation~\eqref{eqf394i094i4095}.
   In other words, the difference is that \m{D_0} is not negligible, and instead
   \Algn{\label{eqh65t4r5t46y7uy6t}D_0 = \TEMPORARYD_{\RHOOO}(D_\DARKCIRC)+ \BIGO{\EPSISMALL p^{-\PADICVALUATION(\ETASMALL!)}}.}
  So upon dividing equation~\eqref{eqj45gj498ug98u} by \m{\TEMPORARYD_{\RHOOO}(D_\DARKCIRC)} and reducing modulo~\m{p} we get that \m{\SCR I} contains
  \Algn{\label{eqi34ir39i093ri0394}
  1\Sqbracketsss{}{KZ}{\FFp}{(\theta^\RHOOO x^{\ETASMALL(p-1)}y^{\rr-\RHOOO(p+1)-\ETASMALL(p-1)}+\ADJUST_{\RHOOO})}.}
  It is easy to show that his represents a generator of \m{\COMPACTIND{KZ}{G}\SUBMODULE(\RHOOO)} (but is trivial in \m{\COMPACTIND{KZ}{G}\QUOTIENTI(\RHOOO)}), which finishes the proof of the desired result as in subsection~\ref{subsecProf1}.

 \subsection{\protect\boldmath If \texorpdfstring{\m{r-\RHOOO(p+1)=1}}{r-\textrho(p+1)=1} and \texorpdfstring{\m{Q}}{Q} is an \texorpdfstring{\m{\infty}}{∞}-dimensional factor of \texorpdfstring{\m{\COMPACTIND{KZ}{G}\QUOTIENTI(\RHOOO)}}{ind(quot(\textrho))}, then \texorpdfstring{\m{Q}}{Q} is a factor of
\texorpdfstring{\m{\COMPACTIND{KZ}{G}\QUOTIENTI(\RHOOO)/\HECKEQ{\RHOOO}}}{ind(quot(\textrho))/T}}
 \label{subsecProf5}
 Let \m{r-\RHOOO(p+1)=1}, so that \Algn{\label{eq6u5jt4iw2ker}\QUOTIENTI(\RHOOO)\cong \sigma_{p-2}(\RHOOO+1).}
 We want to show that \m{\SCR I} contains a representative of a generator of \Algn{\label{eq567ty8i9k}\HECKEQ{\RHOOO}\left(\COMPACTIND{KZ}{G}\QUOTIENTI(\RHOOO)\right).}
 Let %
 \Algn{\label{eqofij4865y439ooij}C_l = \CONS{\RHOOO}{\RHOOO}{l}\Binom{r}{\RHOOO-l}} for \m{l \in \{0,\ldots,\RHOOO\}}. 
 As in subsection~\ref{subsecProf2} we can conclude that
   \Algn{\label{eq9j4u5984u59g8u45t98}\nonumber
      \textstyle & \textstyle
        \nonumber
      \sum_{i(p-1)+\RHOOO\in[0,\RHOOO]} X_i \begin{Lmatrix}
p & 0 \\ 0 & 1
\end{Lmatrix}  \Sqbracketsss{}{KZ}{\QQp}{\ADJUST_{i(p-1)+\RHOOO}}
\\ & \textstyle\null\qquad +  \sum_{i(p-1)+\RHOOO\in(\rr-r+\RHOOO p,\rr]} X_i^\ast\begin{Lmatrix}
1 & 0 \\ 0 & p
\end{Lmatrix} \Sqbracketsss{}{KZ}{\QQp}{\ADJUST^\ast_{\rr-i(p-1)-\RHOOO}}
\nonumber\\ & \textstyle \null\qquad\qquad\phantom{+} {\EQUIVZEROIDEAL}
     \frac{1}{p-1} \sum_{l=0}^{\RHOOO} C_{l}p^l\sum_{\MU\in\FF_p^\times} [\MU]^{-l}\begin{Lmatrix}
p & [\MU] \\ 0 & 1
\end{Lmatrix}  \Sqbracketsss{}{KZ}{\QQp}{\ADJUST_{\RHOOO-l}} \tempheightaer
   + \BIGO{p^{2\ERRR}},
  }
  where
  \Algn{\label{eqfj98fu498ug478gy}\nonumber X_i & \textstyle = p^{-i(p-1)}\Binom{\rr}{i(p-1)+\RHOOO} \Binom{\RHOOO-i}{\RHOOO},\\ X_i^\ast & \textstyle = p^{(i-\RHOOO)(p-1)-1}\Binom{\rr}{i(p-1)+\RHOOO} \Binom{\RHOOO-i}{\RHOOO} .}
  Again, by lemma~\ref{lemmaC1ap},
  \Algn{\label{eqg48u5u849ut54t4554}\PADICVALUATION(X_0) = \PADICVALUATION\left(\Binom{t}{\RHOOO}\right) < \PADICVALUATION\left(p^{-i(p-1)}\Binom{\rr}{i(p-1)+\RHOOO} \Binom{\RHOOO-i}{\RHOOO}\right)=\PADICVALUATION(X_i)}
  for all \m{i} such that \m{i(p-1)+\RHOOO\in[0,\RHOOO)}. By lemma~\ref{lemmaC1bp}, 
  \Algn{\label{eqpodpk3098j293} \PADICVALUATION(X_0) =\PADICVALUATION\left(\Binom{t}{\RHOOO}\right) <\PADICVALUATION\left(p^{(i-\RHOOO)(p-1)-1}\Binom{\rr}{i(p-1)+\RHOOO} \Binom{\RHOOO-i}{\RHOOO}\right) = \PADICVALUATION(X_i^\ast)}
  for all \m{i} such that \m{i(p-1)+\RHOOO\in(\rr-r+\RHOOO p,\rr]}. By lemma~\ref{lemmaC2ap}, 
  \Algn{\label{eqh45875y8t74u849u} \PADICVALUATION(X_0) =\PADICVALUATION\left(\Binom{t}{\RHOOO}\right) <\PADICVALUATION\left(\CONS{\RHOOO}{\RHOOO}{l}\Binom{\rr}{\RHOOO-l}p^l\right)=\PADICVALUATION\left(C_lp^l\right)}
  for all \m{l \in \{1,\ldots,\RHOOO\}}. And, by lemma~\ref{lemmaD},
  \Algn{\label{equgh4j5gi8hghf84ufh} \PADICVALUATION(X_0) =\PADICVALUATION\left(\Binom{t}{\RHOOO}\right)=\PADICVALUATION\left(\Binom{r}{\RHOOO}\right) < 2\ERRR.}
 This means that if we divide both sides of equation~\eqref{eq9j4u5984u59g8u45t98} by \m{\Binom{t}{\RHOOO}} and reduce the resulting  integral element modulo~\m{p}, the two contributing terms are the ``\m{i=0}'' term in the first line of equation~\eqref{eq9j4u5984u59g8u45t98} and the ``\m{l=0}'' term in the third line of equation~\eqref{eq9j4u5984u59g8u45t98}. Therefore \m{\SCR I} contains \Algn{\label{eq409i09wi309ir} \sum_{\MU\in\FF_p} \begin{Lmatrix}
p & [\MU] \\ 0 & 1
\end{Lmatrix}  \Sqbracketsss{}{KZ}{\FFp}{\ADJUST_{\RHOOO}},}
 which is a representative of a generator of \Algn{\label{eq0934ir09it980i}\HECKEQ{\RHOOO}\left(\COMPACTIND{KZ}{G}\QUOTIENTI(\RHOOO)\right),} and that completes the proof.

\begin{proof}{Proof of theorem~\ref{MAINTHEOREMp} (\m{\Leftrightarrow} theorem~\ref{MAINTHEOREM}). }
Let \m{Q} be an infinite-dimensional factor of \m{\overline{\Theta}_{\rr+2,\aaa}}.
Subsections~\ref{subsecProf1} and~\ref{subsecProf2} imply the following two facts about \m{Q}.

\emph{1.} \m{Q} is not a factor of {\m{\FILT N_\alpha}} for {\m{\alpha\in\{0,\ldots,\RHOOO-1\}}}.

\emph{2.}
\m{Q} is not a factor of {\m{\FILT N_\alpha}}  for
\m{\alpha\in\{\RHOOO+1,\ldots,\TEMPORARYI-1\}}.

From these two facts we can conclude that either \m{Q} is a factor of \m{\COMPACTIND{KZ}{G}\SUBMODULE(\RHOOO)} or it is a factor of \m{\COMPACTIND{KZ}{G}\QUOTIENTI(\RHOOO)}. In light of theorem~\ref{BERGERX1} and as in section~10 of~\cite{Ars1} this implies equation~\eqref{qowur983u4t3yt8374t4y}
for \m{r-\RHOOO(p+1) \in \{-1,\ldots,p-1\}\backslash\{-1,1,p-2\}}.

Subsections~\ref{subsecProf3} and~\ref{subsecProf5} prove the following facts for \m{r-\RHOOO(p+1)\in\{1,p-2\}}.

\emph{3.} 
If \Algn{\label{eqwe-or234r983uit9845}r-\RHOOO(p+1)=p-2} (and therefore \m{\SUBMODULE(\RHOOO)\cong \sigma_{p-2}(\RHOOO)}) then \m{Q} is a factor of \m{\COMPACTIND{KZ}{G}\QUOTIENTI(\RHOOO)}.

\emph{4.} 
If \Algn{\label{eqw0ire98rui983u45g98u45}r-\RHOOO(p+1)=1} (and therefore \m{\QUOTIENTI(\RHOOO)\cong \sigma_{p-2}(\RHOOO+1)}) and \m{Q} is a factor of \m{\COMPACTIND{KZ}{G}\QUOTIENTI(\RHOOO)} then \m{Q} is a factor of of \m{\COMPACTIND{KZ}{G}\QUOTIENTI(\RHOOO)/\HECKEQ{\RHOOO}}.

These two claims imply equation~\eqref{qowur983u4t3yt8374t4y} for \m{r-\RHOOO(p+1)\in\{1,p-2\}}. Since we assume that \m{p+1\nmid k-1}, we have \m{r-\RHOOO(p+1)\ne-1}, and therefore the proof is complete.
\end{proof}

 \begin{proof}{Proof that theorem~\ref{MAINTHEOREM} implies corollary~\ref{MAINCORO}. }
As in subsection~4.2 of~\cite{BLZ04} we can use theorem~\ref{MAINTHEOREM} to conclude that \m{\mu_l} is supported on \Algn{\label{eq43u98ur4398r}\textstyle \left[0,\frac{1}{p+1}+\frac{\log_pl}{l-1}\right]\cup \left[\frac{p}{p+1}-\frac{\log_pl}{l-1},1\right],} and that completes the proof because \Algn{\label{eq340j4gj948g5i}\frac{\log_pl}{l-1}\to 0} as \m{l\to\infty}---we omit the details.
\end{proof}

\section{Combinatorics}
\label{secCom}

\begin{lemma}\label{lemmaD}
If \m{\alpha \in \ZZ_{\GEQSLANT1}} and \m{\beta \in \{0,\ldots,\alpha\}} then
\Algn{\label{eq9043i0r} \PADICVALUATION\left(\Binom{\alpha}{\beta}\right) \LEQSLANT \lfloor\log_p \alpha\rfloor.}
\end{lemma}

\begin{proof}{Proof. }
A theorem by Kummer says that \Algn{\label{qwe09i3908ri3894ut3ut}\PADICVALUATION\left(\Binom{\alpha}{\beta}\right)} is the number of times one carries over a digit when adding \m{\beta} and \m{\alpha-\beta}, and is therefore strictly less than the number \m{\lfloor \log_p \alpha\rfloor+1} of digits of \m{\alpha}. %
\end{proof}

Let \m{\alpha\in\{\RHOOO,\ldots,\TEMPORARYI-1\}},  let
\Algn{\label{eqr092ir09ir903} \RHOOOPRIME=\lceil\frac{r-\alpha}{p}\rceil-1,}
and for \m{l\in\{\alpha-\RHOOOPRIME,\ldots,\alpha\}} let
\Algn{\label{eqori230ri349ut9348i34} C_l = \CONS{\RHOOOPRIME}{\alpha}{l}\Binom{r}{\alpha-l}.}
Note that \m{\RHOOOPRIME\LEQSLANT\RHOOO} and if \m{r=\RHOOOPRIME(p+1)+1} and  \m{\alpha=\RHOOO} then \m{\RHOOOPRIME=\RHOOO}.
   The constants \m{C_l} are precisely those constants that satisfy
   \Algn{\label{eqr34t4yu67i67i76i6} \sum_{l=\alpha-\RHOOOPRIME}^\alpha C_l \Binom{r}{\alpha-l}^{-1}\Binom{(p-1)X+\alpha}{\alpha-l} = \Binom{\RHOOOPRIME-X}{\RHOOOPRIME} \in \QQ_p[X].}
   Moreover, let \Algn{\label{eq65u67i78oi78o78i78i}\nonumber X_i & \textstyle = p^{-i(p-1)}\Binom{r}{i(p-1)+\alpha} \Binom{\RHOOOPRIME-i}{\RHOOOPRIME} \text{ for } i \in\ZZ%
    ,\\ X_i^\ast & \textstyle = p^{i(p-1)+2\alpha-r}\Binom{r}{i(p-1)+\alpha} \Binom{\RHOOOPRIME-i}{\RHOOOPRIME} \text{ for } i \in\ZZ%
    .}

\begin{lemma}
\label{lemmaC1a}
If \m{\alpha>\RHOOO} and \m{i(p-1)+\alpha\in[0,\alpha)} then 
\Algn{\label{eq0fei09fi904i} \PADICVALUATION\left(X_0\right) < \PADICVALUATION\left(X_i\right).}
\end{lemma}

\begin{proof}{Proof. }
Note that \m{i<0}, so let us write \m{j=-i>0}. We have
\Algn{\label{eq09ri093ir9034} X_0 & \nonumber\textstyle = \binom{r}{\alpha},\\ \textstyle X_i & \textstyle= p^{j(p-1)}\Binom{r}{\alpha-j(p-1)} \Binom{\RHOOOPRIME+j}{\RHOOOPRIME} = p^{j(p-1)}\Binom{r}{\alpha} \frac{\alpha_{j(p-1)}}{(r-\alpha+j(p-1))_{j(p-1)}} \Binom{\RHOOOPRIME+j}{j}.}
Therefore we want to show that
\Algn{\label{eq309ri3209rit3904ti} \PADICVALUATION\left((r-\alpha+j(p-1))_{j(p-1)}\right) < \PADICVALUATION\left(\alpha_{j(p-1)}\right) + \PADICVALUATION\left((\RHOOOPRIME+j)_j\right) + j(p-1)-\PADICVALUATION(j!).}
Note that \m{\PADICVALUATION(j!) \LEQSLANT \frac{j}{p-1}}, so it is enough to show that
\Algn{\label{eq3ru3498tu59348ut948ytu} \PADICVALUATION\left((r-\alpha+j(p-1))_{j(p-1)}\right) < \PADICVALUATION\left(\alpha_{j(p-1)}\right) + \PADICVALUATION\left((\RHOOOPRIME+j)_j\right) + j\left(p-1-\frac{1}{p-1}\right).}
We have
\Algn{\label{eqwepr923ir093i90tr} \lfloor\frac{r-\alpha+j(p-1)}{p}\rfloor \LEQSLANT \RHOOOPRIME+j = \lceil\frac{r-\alpha+jp}{p}\rceil-1.}
We also have
\Algn{\label{eqr9i2039rti390it} \lceil\frac{r-\alpha+1}{p}\rceil \GEQSLANT \RHOOOPRIME+1 = \lceil\frac{r-\alpha}{p}\rceil.}
Therefore each term of the product
\Algn{\label{eqre-02ir093i4t90i} (r-\alpha+j(p-1))_{j(p-1)} = (r-\alpha+j(p-1)) \cdots (r-\alpha+1)} that is divisible by \m{p} is \m{p} times a term of the product
\Algn{\label{eqwe902ir9034t0}(\RHOOOPRIME+j)_j=(\RHOOOPRIME+j)\cdots(\RHOOOPRIME+1),}
implying that
\Algn{\label{eqr23r34t45y56u56uu6} \PADICVALUATION\left((r-\alpha+j(p-1))_{j(p-1)}\right) \LEQSLANT \PADICVALUATION\left((\RHOOOPRIME+j)_j\right) +w,}
where \m{w} is the number of terms of the product in equation~\eqref{eqre-02ir093i4t90i} that are divisible by \m{p}. Equation~\eqref{eq3ru3498tu59348ut948ytu} follows from the fact that
\Algn{\label{eqr-09i3r09i34t9u43598tyu459yu54y} w\LEQSLANT \lfloor\frac{j(p-1)+p-1}{p}\rfloor < j \left(p-1-\frac{1}{p-1}\right).}
\end{proof}

\begin{lemma}
\label{lemmaC1b}
If \m{\alpha>\RHOOO} and \m{i(p-1)+\alpha\in(\RHOOOPRIME(p-1)+\alpha,r]} then 
\Algn{\label{eq23r34t45y54y56u65} \PADICVALUATION\left(X_0\right) < \PADICVALUATION\left(X_i^\ast\right).}
\end{lemma}

\begin{proof}{Proof. }
We have \m{i(p-1)+\alpha\in(\RHOOOPRIME(p-1)+\alpha,r]} and therefore \Algn{\label{ewqt4t45ty45y56u56u}j=i(p-1)+2\alpha-r\in \left[\lceil\frac{r-\alpha}{p}\rceil(p-1)+2\alpha-r,\alpha\right]\subseteq \left[\frac{(p+1)\alpha-r}{p},\alpha\right] \subseteq (0,\alpha].} So we can write
\Algn{\label{eq34y56u67u76i7i7i} X_0 & \nonumber\textstyle = \binom{r}{\alpha},\\ \textstyle X_i^\ast & \textstyle= p^{i(p-1)+2\alpha-r}\Binom{r}{i(p-1)+\alpha} \Binom{\RHOOOPRIME-i}{\RHOOOPRIME} = p^{j}\Binom{r}{\alpha} \frac{\alpha_{j}}{(r-\alpha+j)_{j}} \Binom{\RHOOOPRIME-i}{\RHOOOPRIME}.}
Therefore we want to show that
\Algn{\label{eqr34t456u67u67u} \PADICVALUATION\left((r-\alpha+j)_{j}\right) < \PADICVALUATION\left(\alpha_{j}\right) + \PADICVALUATION\left(\binom{\RHOOOPRIME-i}{\RHOOOPRIME}\right) + j.}
We have
\Algn{\label{eqwe09i2390ri309tui3498t} \binom{\RHOOOPRIME-i}{\RHOOOPRIME} = (-1)^{\RHOOOPRIME}\binom{i-1}{\RHOOOPRIME}=(-1)^{\RHOOOPRIME}\binom{i-1}{i-\RHOOOPRIME-1}=(-1)^{\RHOOOPRIME}\frac{(i-1)\cdots(\RHOOOPRIME+1)}{(i-\RHOOOPRIME-1)!}.}
Because \m{\alpha>\RHOOOPRIME} and \m{i>\RHOOOPRIME}, we have \m{r\LEQSLANT ip+\alpha-1}, and therefore \m{\alpha-j\LEQSLANT i-1}. Moreover, \m{\alpha\GEQSLANT \RHOOO+1\GEQSLANT\RHOOOPRIME+1}. Therefore the union of the intervals \m{(\RHOOOPRIME+1,i-1]} and \m{(\alpha-j+1,\alpha]} is a single interval. This implies that \Algn{\label{eqw09ri90i3409ti} \PADICVALUATION\left(\alpha_j (i-1)_{i-\RHOOOPRIME-1}\right)\GEQSLANT \PADICVALUATION\left(\max\{i-1,\alpha\}\cdots \min\{\RHOOOPRIME+1,\alpha-j\}\right).}
We also have
\Algn{\label{eqw9turi349tu89t}\PADICVALUATION\left((i-\RHOOOPRIME-1)!\right)\LEQSLANT \frac{i-\RHOOOPRIME-1}{p-1}.}
Since \Algn{\label{eqwei09ti3409tit}p\max\{i-1,\alpha\}+p-1\GEQSLANT r-\alpha+j \text{ and }p\min\{\RHOOOPRIME+1,\alpha-j\}-p+1\LEQSLANT r-\alpha+1,} we have
\Algn{\label{eqw09i23r0934ir903tui98}\PADICVALUATION\left((r-\alpha+j)_{j}\right) \LEQSLANT \PADICVALUATION\left(\max\{i-1,\alpha\}\cdots \min\{\RHOOOPRIME+1,\alpha-j\}\right) + \frac{j+p-1}{p}.}
So it is enough to show that
\Algn{\label{eqwe0-2or093it0943} j > \frac{j+p-1}{p} + \frac{i-\RHOOOPRIME-1}{p-1} ,} which follows from \Algn{\label{eqwe09i2r9i349it} \frac{(i(p-1)+2\alpha-r)(p-1)}{p} \GEQSLANT \frac{i-\RHOOOPRIME+p-2}{p-1}.} The latter inequality follows from the fact that \m{\gamma(i)=\frac{(i(p-1)+2\alpha-r)(p-1)}{p}-\frac{i-\RHOOOPRIME+p-2}{p-1}} is increasing in \m{i} and \m{\gamma(\RHOOOPRIME+1) \GEQSLANT \frac{2(p-1)}{p}-1>0}.
\end{proof}

\begin{lemma}
\label{lemmaC2a}
If \m{\alpha>\RHOOO} and \m{l\in\{\alpha-\RHOOOPRIME,\ldots,\alpha\}} then
\Algn{\label{eqiofj5jg498j5}
\PADICVALUATION\left(X_0\right) < \PADICVALUATION\left(C_lp^l\right).} 
\end{lemma}

\begin{proof}{Proof. }
For \m{l\in\{\alpha-\RHOOOPRIME,\ldots,\alpha\}} let \Algn{\label{eq398rtu349tu3498tu45}C_l^\prime = C_l \Binom{r}{\alpha-l}^{-1}\in\ZZ_p,} so that
\Algn{\label{eq56j678ik78o89l89o} \sum_{j=\alpha-\RHOOOPRIME}^\alpha C_j^\prime \Binom{(p-1)X+\alpha}{\alpha-j} = \Binom{\RHOOOPRIME-X}{\RHOOOPRIME} \in \QQ_p[X].}
We have
\Algn{\label{eqgtryjukilo98i7u6y5t4} X_0 & \nonumber\textstyle = \binom{r}{\alpha},\\ \textstyle C_lp^l & \textstyle= C_l^\prime \binom{r}{\alpha-l}p^l = C_l^\prime \binom{r}{\alpha}\frac{\alpha_{l}}{(r-\alpha+l)_l}p^l.}
So we want to show that
\Algn{\label{eqwe-i209i4309ti390ti4} \PADICVALUATION\left(C_l^\prime\right) > \PADICVALUATION\left((r-\alpha+l)_l\right)-\PADICVALUATION\left(\alpha_l\right)-l.}
First suppose that \m{\alpha(p+1)+p-1\GEQSLANT r+l}. Then \m{p\alpha+p-1\GEQSLANT r-\alpha+l}, implying that the largest term in the sequence \m{(r-\alpha+l,\ldots,r-\alpha+1)} that is divisible by \m{p} is at most as large as \m{p} times the largest term in the sequence \m{(\alpha,\ldots,\alpha-l+1)}. 
Moreover, \Algn{\label{eqwe098u2r98u389ut934} p\alpha -pl +1\LEQSLANT p\alpha-p(\alpha-\RHOOOPRIME)+1=p\RHOOOPRIME+1\LEQSLANT r-\alpha+1,}
implying that the smallest term in the sequence \m{(r-\alpha+l,\ldots,r-\alpha+1)} that is divisible by \m{p} is at least as large as \m{p} times the smallest term in the sequence \m{(\alpha,\ldots,\alpha-l+1)}. Therefore
\Algn{\label{eeqwe09i349t349ti3498} \PADICVALUATION\left((r-\alpha+l)_l\right)-\PADICVALUATION\left(\alpha_l\right)-l\LEQSLANT w-l \LEQSLANT \lfloor\frac{l+p-1}{p}\rfloor -l\LEQSLANT 0,}
where \m{w} is the number of terms in the sequence \m{(r-\alpha+l,\ldots,r-\alpha+1)} that are divisible by \m{p}.
Equality in equation~\eqref{eeqwe09i349t349ti3498} holds if and only if \m{l=1} and \m{r-\alpha+1=p\alpha}, which can never happen. 
So the right side of equation~\eqref{eqwe-i209i4309ti390ti4} is negative, implying equation~\eqref{eqiofj5jg498j5} in the case when \m{\alpha(p+1)+p-1\GEQSLANT r+l}.
Now suppose that \m{\alpha(p+1)+p\LEQSLANT r+l}, so that \Algn{\label{eqwe-09i34rt0934i409it}l \GEQSLANT p(\alpha-\RHOOOPRIME)+p-1 \GEQSLANT 2p-1.} For \m{j\in\{\alpha-\RHOOOPRIME,\ldots,\alpha\}} let \Algn{\label{eqwe09i23r093i3904tiu948t}C_j^\mydprime = (-1)^{\RHOOOPRIME} (p-1)^{j-\alpha} C_j^\prime \RHOOOPRIME_{j-\alpha+\RHOOOPRIME} \in\ZZ_p,} so that
\Algn{\label{eq45t56u78ik78i78} \sum_{j=\alpha-\RHOOOPRIME}^\alpha C_j^\mydprime \prod_{u=0}^{\alpha-j-1} \left(X+\frac{\alpha-u}{p-1}\right) = (X-1)\cdots(X-\RHOOOPRIME) \in \QQ_p[X].}
We use the fact that among any \m{l} consecutive integers there can be at most one whose valuation is at least \m{\log_p l}. Moreover, if there is such a term then the sum of the valuations of all the other terms is at most \m{\frac{l-1}{p-1}}, and if there is no such term then the sum of the valuations of all the terms is at most \m{\frac{l-1}{p-1}+\lfloor \log_pl\rfloor}. So if there is no term in the sequence \m{(r-\alpha+l,\ldots,r-\alpha+1)} 
whose valuation is at least \m{\log_p l}, then 
\Algn{\label{eq0954it490it4598i90ti} \PADICVALUATION\left((r-\alpha+l)_l\right)-l < \lfloor \log_pl\rfloor - \frac{l(p-2)}{p-1}  \LEQSLANT 0,}
so \Algn{\label{eq390i40934i3094tir} \PADICVALUATION\left(C_l^\prime \binom{r}{\alpha}\frac{\alpha_{l}}{(r-\alpha+l)_l}p^l\right) \GEQSLANT \PADICVALUATION\left( \binom{r}{\alpha}\frac{1}{(r-\alpha+l)_l}p^l\right) > \PADICVALUATION\left( \binom{r}{\alpha}\right).} Suppose now that there is a term in the sequence \m{(r-\alpha+l,\ldots,r-\alpha+1)} 
whose valuation is \m{\gamma\GEQSLANT \log_p l}. If \m{\gamma \LEQSLANT \frac{l(p-2)}{p-1}} then we can similarly deduce that
\Algn{\label{eq3984tu498tu4598tu4598ut4598y} \PADICVALUATION\left(C_l^\prime \binom{r}{\alpha}\frac{\alpha_{l}}{(r-\alpha+l)_l}p^l\right) \GEQSLANT \PADICVALUATION\left( \binom{r}{\alpha}\frac{1}{(r-\alpha+l)_l}p^l\right) > \PADICVALUATION\left( \binom{r}{\alpha}\right),}
so suppose that \m{\gamma \GEQSLANT \frac{l(p-2)}{p-1}\GEQSLANT \log_p l.} If the term in  \m{(r-\alpha+l,\ldots,r-\alpha+1)} whose valuation is \m{\gamma} is \m{p} times a term in \m{(\alpha,\ldots,\alpha-l+1)}, then we can similarly deduce that
\Algn{\label{eqr34t456y56u67u67u} \PADICVALUATION\left(C_l^\prime \binom{r}{\alpha}\frac{\alpha_{l}}{(r-\alpha+l)_l}p^l\right) \GEQSLANT \PADICVALUATION\left( \binom{r}{\alpha}\frac{\alpha_l}{(r-\alpha+l)_l}p^l\right) > \PADICVALUATION\left( \binom{r}{\alpha}\right),}
so suppose that the term in \m{(r-\alpha+l,\ldots,r-\alpha+1)} whose valuation is \m{\gamma} is in the subsequence \m{(r-\alpha+l,\ldots,p(\alpha+1))}. Let \m{q=p^{\gamma-1}}. Informally speaking, the assumptions that \m{q} is a ``large'' power of \m{p} and that  \m{(p(\RHOOOPRIME+1)+l,\ldots,p(\alpha+1))} contains a multiple of \m{pq} (say \m{bpq}) and that \m{\alpha>\RHOOOPRIME} imply that \m{\alpha} and \m{\RHOOOPRIME} are ``just below'' a multiple of \m{q}, and \Algn{\label{eqwe2r983ur9384u983tu} \RHOOOPRIME+l/p +1\GEQSLANT bq > \alpha > \RHOOOPRIME \text{ and } q\GEQSLANT p^{\frac{l(p-2)}{p-1}-1} \text{ and } l\GEQSLANT 2p-1.}
For \m{u\in\{0,\ldots,\RHOOOPRIME-1\}} let \m{z_u} be the integer in \m{\{1,\ldots,q\}} that is congruent to \m{\frac{u-\alpha}{p-1}} modulo \m{q}, and 
for \m{j\in\{\alpha-\RHOOOPRIME,\ldots,\alpha\}} let \m{C_j^\mytrprime} be such that
\Algn{\label{eq45t456u67u6i67i} \sum_{j=\alpha-\RHOOOPRIME}^\alpha C_j^\mytrprime \prod_{u=0}^{\alpha-j-1} \left(X-z_u\right) = \prod_{i=1}^{q-s}(X-i)^b \prod_{i=q-s+1}^{q} (X-i)^{b-1} \in \QQ_p[X],} 
where \m{s=bq-\RHOOOPRIME}.
By reducing equation~\eqref{eq45t56u78ik78i78}  modulo \m{q} we get \m{C_j^\mydprime \equiv C_j^\mytrprime \bmod q}. 
For \m{j \in \{\alpha-\RHOOOPRIME,\ldots,\alpha\}}, let \m{F_j(X)=\prod_{u=0}^{\alpha-j-1} \left(X-z_u\right)}. Then 
\Algn{\label{eq98ur93utu8349t}F_\alpha(X)\mid \cdots\mid F_{\alpha-\RHOOOPRIME}(X),} and if \m{i_0\in\{\alpha-\RHOOOPRIME,\ldots,\alpha\}} is the smallest index such that \m{F_{i_0}(X)} divides \Algn{\label{eqwe092ir09ir3904t}G(X)=\prod_{i=1}^{q-s}(X-i)^b \prod_{i=q-s+1}^{q} (X-i)^{b-1},} then equation~\eqref{eq45t456u67u6i67i} implies that \m{C_j^\mytrprime=0} for all \m{j \in \{i_0+1,\ldots,\alpha\}}, and therefore \m{C_j^\mydprime\equiv 0 \bmod q} for all \m{j \in \{i_0+1,\ldots,\alpha\}}. We want to show that \m{i_0>l}, \IE that
\Algn{\label{eqwe092ir092ir093i3409t} \prod_{u=0}^{\alpha-l} \left(X-z_u\right) \mid G(X).}
We have
\Algn{\label{eqfergh56j67i7i78i78i} \prod_{u=0}^{(b-1)q-1} \left(X-z_u\right) = \prod_{i=1}^{q}(X-i)^{b-1},}
so in order to show equation~\eqref{eqwe092ir092ir093i3409t} it is enough to show that
\Algn{\label{eqwe9i23r09i230ri3094ti}\prod_{u=(b-1)q}^{\alpha-l} \left(X-z_u\right)\mid \prod_{i=1}^{q-s}(X-i),}
\IE that the sets
\Algn{\label{eqwe-023ir09ri309ri} \left\{bq-\alpha+i \,|\, i \in \{0,\ldots,\alpha-l-(b-1)q\}\right\}} and \Algn{\label{eqwe092i093i49tu3945u8954} \left\{q-(p-1)i \,|\, i \in \{0,\ldots,bq-\RHOOOPRIME-1\}\right\}} are disjoint. This follows from
\Algn{\label{eqwe092ir3904ti349tu}\nonumber
& \textstyle l \GEQSLANT p(bq-\RHOOOPRIME-1) \\ &\nonumber \null\qquad\textstyle \Longrightarrow l > (p-1)(bq-\RHOOOPRIME-1)
\\ & \null\qquad\textstyle \Longrightarrow q-l < q-(p-1)(bq-\RHOOOPRIME-1).} So \m{i_0 > l}, and therefore \m{C_l^\mydprime\equiv 0 \bmod q}. Since
\Algn{\label{eqr345y678o990[90[0[0-} X_0 & \nonumber\textstyle = \binom{r}{\alpha},\\ \textstyle C_lp^l & \textstyle= C_l^\prime \binom{r}{\alpha}\frac{\alpha_{l}}{(r-\alpha+l)_l}p^l = C_l^\mydprime \binom{r}{\alpha}\frac{\alpha\cdots (\RHOOOPRIME+1)}{(r-\alpha+l)_l}p^l,}
and since
\m{\PADICVALUATION(C_l^\mydprime)\GEQSLANT \gamma-1}  and 
\Algn{\label{eqwe09i2309ri290}\PADICVALUATION\left((r-\alpha+l)_l\right) \LEQSLANT \gamma+\frac{l-1}{p-1}}
and
\Algn{\label{eqwe902ir09ri9034ti} l-\frac{l-1}{p-1}-1 = \frac{(l-1)(p-2)}{p-1}\GEQSLANT 2(p-2)>0,}
we can deduce equation~\eqref{eqiofj5jg498j5} and complete the proof.
\end{proof}

\begin{lemma}
\label{lemmaC1ap}
If \m{r=\RHOOO(p+1)+1} and \m{\alpha=\RHOOO} and \m{i(p-1)+\RHOOO\in[0,\RHOOO)} then 
\Algn{\label{eqofij4tytgujwis} \PADICVALUATION\left(X_0\right) < \PADICVALUATION\left(X_i\right).}
\end{lemma}

\begin{proof}{Proof. }
The proof is similar to the proof of lemma~\ref{lemmaC1a}. Note that \m{i<0}, so let us write \m{j=-i>0}. We have
\Algn{\label{eq092ir39ut98uty8945yu8457y} X_0 & \nonumber\textstyle = \binom{r}{\RHOOO},\\ \textstyle X_i & \textstyle= p^{j(p-1)}\Binom{r}{\RHOOO-j(p-1)} \Binom{\RHOOO+j}{\RHOOO} = p^{j(p-1)}\Binom{r}{\RHOOO} \frac{\RHOOO_{j(p-1)}}{(\RHOOO p+j(p-1)+1)_{j(p-1)}} \Binom{\RHOOO+j}{j}.}
Therefore we want to show that
\Algn{\label{eqrwty65u778o89p} \PADICVALUATION\left((\RHOOO p+j(p-1)+1)_{j(p-1)}\right) < \PADICVALUATION\left((\RHOOO+j)_{jp}\right) + j\left(p-1-\frac{1}{p-1}\right),}
since \m{\PADICVALUATION(j!) \LEQSLANT \frac{j}{p-1}}.
Since \Algn{\label{eq9q38ru394ur3894t}(\RHOOO+j)p > \RHOOO p+j(p-1)+1 \text{ and }\RHOOO p+2>(\RHOOO-j(p-1)+1)p,}
 each term of the product
\Algn{\label{eq4t34t45y4y45t6u67u} (\RHOOO p+j(p-1)+1)_{j(p-1)}} that is divisible by \m{p} is \m{p} times a term of the product
\Algn{\label{eqrwtret45y45y}(\RHOOO+j)_{jp}.}
This together with
\Algn{\label{eqwe9i23ri3209tri390t4i} j\left(p-1-\frac{1}{p-1}\right) > \frac{(j+1)(p-1)}{p}} implies equation~\eqref{eqrwty65u778o89p}.
\end{proof}

\begin{lemma}
\label{lemmaC1bp}
If \m{r=\RHOOO(p+1)+1} and \m{\alpha=\RHOOO} and \m{i(p-1)+\RHOOO\in(\RHOOO p,r]} then
\Algn{\label{eqrritbgbvfdsver} \PADICVALUATION\left(X_0\right) < \PADICVALUATION\left(X_i^\ast\right).}
\end{lemma}

\begin{proof}{Proof. } The proof is similar to the proof of lemma~\ref{lemmaC1b}. 
We have \m{i(p-1)+\RHOOO\in(\RHOOO p,r]} and therefore \Algn{\label{eqk8t78io7ti7to8}j=(i-\RHOOO)(p-1)-1\in \left[ p-2,\RHOOO\right].} So we can write
\Algn{\label{eq78o7o78o89o8o8o8p} X_0 & \nonumber\textstyle = \binom{r}{\RHOOO},\\ \textstyle X_i^\ast & \textstyle=  p^{j}\Binom{r}{\RHOOO} \frac{\RHOOO_{j}}{(\RHOOO p+j+1)_{j}} \Binom{\RHOOO-i}{\RHOOO}=  (-1)^\RHOOO p^{j}\Binom{r}{\RHOOO} \frac{(i-1)_{i-\RHOOO+j-1}}{(\RHOOO p+j+1)_{j} (i-\RHOOO-1)!},} where the last equality follows as in the proof of lemma~\ref{lemmaC1b}. 
Therefore we want to show that
\Algn{\label{equ67u765it78it78i} \PADICVALUATION\left((\RHOOO p+j+1)_{j}\right) + \frac{i-\RHOOO-1}{p-1} < \PADICVALUATION\left((i-1)_{i-\RHOOO+j-1}\right) + j.}
Since \m{ip-1\GEQSLANT i(p-1)+\RHOOO} and \m{\RHOOO p - jp+1\LEQSLANT \RHOOO p+2}, we have
\Algn{\label{eq09it3409it3904t}\PADICVALUATION\left((\RHOOO p+j+1)_{j}\right) \LEQSLANT \PADICVALUATION\left((i-1)_{i-\RHOOO+j-1}\right) + \frac{j+p-1}{p},}
and equation~\eqref{equ67u765it78it78i} follows from
\Algn{\label{eqwr09i3490u39tu485ut} j > \frac{j+p-1}{p} + \frac{i-\RHOOO-1}{p-1},} as in 
the proof of lemma~\ref{lemmaC1b}.
\end{proof}

\begin{lemma}
\label{lemmaC2ap}
If \m{r=\RHOOO(p+1)+1} and \m{\alpha=\RHOOO} and \m{l\in\{1,\ldots,\RHOOO\}} then 
\Algn{\label{eqrjenigtr5y48ujfweiveineriu}
\PADICVALUATION\left(X_0\right) < \PADICVALUATION\left(C_lp^l\right).} 
\end{lemma}

\begin{proof}{Proof. }
For \m{l\in\{0,\ldots,\RHOOO\}} let \Algn{\label{eq5h64756j57jik}C_l^\prime = C_l \Binom{r}{\RHOOO-l}^{-1}\in\ZZ_p,} so that
\Algn{\label{eq45y56u56j57jk678k} \sum_{j=0}^\RHOOO C_j^\prime \Binom{(p-1)X+\RHOOO}{\RHOOO-j} = \Binom{\RHOOO-X}{\RHOOO} \in \QQ_p[X].}
We have
\Algn{\label{eq46j578k678k8k79l} X_0 & \nonumber\textstyle = \binom{r}{\RHOOO},\\ \textstyle C_lp^l & \textstyle= C_l^\prime \binom{r}{\RHOOO-l}p^l = C_l^\prime \binom{r}{\RHOOO}\frac{\RHOOO_{l}}{(\RHOOO p+l+1)_l}p^l.}
Therefore in order to prove the lemma it is enough to show that, for \m{l\in\{1,\ldots,\RHOOO\}},
\Algn{\label{eq4ye56u56j7k68k879} \PADICVALUATION\left(C_l^\prime\right) > \PADICVALUATION\left((\RHOOO p+l+1)_l\right)-\PADICVALUATION\left(\RHOOO_l\right)-l.}
Note that equation~\eqref{eq4ye56u56j7k68k879} is not true for \m{l=0}, since both sides are zero. As in the proof of lemma~\ref{lemmaC2a}, if \m{l\LEQSLANT p-2} then we can deduce that
\Algn{\label{eqi0934ti309ti4509} \PADICVALUATION\left((\RHOOO p+l+1)_l\right)-\PADICVALUATION\left(\RHOOO_l\right)-l \LEQSLANT \frac{l+p-1}{p}-l\LEQSLANT 0,}
with equality if and only if \m{l=1} and \m{\RHOOO p+2=\RHOOO p}, which can never happen. So let us assume that \m{l\GEQSLANT p-1}. 
For \m{j\in\{0,\ldots,\RHOOO\}} let \Algn{\label{eq40t9i40ti4095yi506k59}C_j^\mydprime = (-1)^{\RHOOO} (p-1)^{j-\RHOOO} C_j^\prime \RHOOO_{j} \in\ZZ_p,} so that
\Algn{\label{eq30kf03ik409gi5ik95hk} \sum_{j=0}^\RHOOO C_j^\mydprime \prod_{u=0}^{\RHOOO-j-1} \left(X+\frac{\RHOOO-u}{p-1}\right) = (X-1)\cdots(X-\RHOOO) \in \QQ_p[X].}
As in the proof of lemma~\ref{lemmaC2a}, if there is no term in \m{(\RHOOO p+l+1,\ldots,\RHOOO p+2)} 
whose valuation is at least \m{\log_p l}, then 
\Algn{\label{eq6i78i6o9l78979o} \PADICVALUATION\left((\RHOOO p+l+1)_l\right)-l < \lfloor \log_pl\rfloor - \frac{l(p-2)}{p-1}  \LEQSLANT 0,}
implying equation~\eqref{eq4ye56u56j7k68k879}.  Suppose now that there is a term in  \Algn{\label{eq3-04ot3-0to3}(\RHOOO p + l + 1,\ldots,\RHOOO p + 2)} 
whose valuation is \m{\gamma\GEQSLANT \log_p l}. If \m{\gamma \LEQSLANT \frac{l(p-2)}{p-1}} then we can similarly deduce equation~\eqref{eq4ye56u56j7k68k879}
as in the proof of lemma~\ref{lemmaC2a},
so suppose that \m{\gamma \GEQSLANT \frac{l(p-2)}{p-1}\GEQSLANT \log_p l.} Let \m{q=p^{\gamma-1}}. We have, for some positive integer \m{b}, \Algn{\label{eq45eo-l5o6hykl560ok} \RHOOO+(l+1)/p \GEQSLANT bq >  \RHOOO \text{ and } q\GEQSLANT p^{\frac{l(p-2)}{p-1}-1} \text{ and } l\GEQSLANT p-1.}
For \m{u\in\{0,\ldots,\RHOOO-1\}} let \m{z_u} be the integer in \m{\{1,\ldots,q\}} that is congruent to \m{\frac{u-\RHOOO}{p-1}} modulo \m{q}, and 
for \m{j\in\{0,\ldots,\RHOOO\}} let \m{C_j^\mytrprime} be such that
\Algn{\label{eq56ujr6ujt67ikt7k7tkt7k} \sum_{j=0}^\RHOOO C_j^\mytrprime \prod_{u=0}^{\RHOOO-j-1} \left(X-z_u\right) = \prod_{i=1}^{q-s}(X-i)^b \prod_{i=q-s+1}^{q} (X-i)^{b-1} \in \QQ_p[X],} 
where \m{s=bq-\RHOOO}.
By reducing equation~\eqref{eq30kf03ik409gi5ik95hk}  modulo \m{q} we get \m{C_j^\mydprime \equiv C_j^\mytrprime \bmod q}. 
For \m{j \in \{0,\ldots,\RHOOO\}}, let \m{F_j(X)=\prod_{u=0}^{\RHOOO-j-1} \left(X-z_u\right)}. Then 
\Algn{\label{eq4yh5e6yr5ur6u7}F_\RHOOO(X)\mid \cdots\mid F_{0}(X),} and if \m{i_0\in\{0,\ldots,\RHOOO\}} is the smallest index such that \m{F_{i_0}(X)} divides \Algn{\label{eq5ur67uj67j67kjtk}G(X)=\prod_{i=1}^{q-s}(X-i)^b \prod_{i=q-s+1}^{q} (X-i)^{b-1},} then equation~\eqref{eq56ujr6ujt67ikt7k7tkt7k} implies that \m{C_j^\mytrprime=0} for all \m{j \in \{i_0+1,\ldots,\RHOOO\}}, and therefore \m{C_j^\mydprime\equiv 0 \bmod q} for all \m{j \in \{i_0+1,\ldots,\RHOOO\}}. So again, as in the proof of lemma~\ref{lemmaC2a},  we can complete the proof by noting that
\Algn{\label{eqe7ity8it78i78oly89ly8}\nonumber
& \textstyle l \GEQSLANT p(bq-\RHOOO)-1 \\ &\nonumber \null\qquad\textstyle \Longrightarrow l > (p-1)(bq-\RHOOO-1)
\\ & \null\qquad\textstyle \Longrightarrow q-l < q-(p-1)(bq-\RHOOO-1)\nonumber
\\ & \null\qquad\textstyle \Longrightarrow i_0>l,}
and therefore that, due to equation~\eqref{eq46j578k678k8k79l}, equation~\eqref{eqrjenigtr5y48ujfweiveineriu} follows from
\Algn{\label{eq094ik309ti09ti} &\nonumber \textstyle\PADICVALUATION\left(C_l^\mydprime\right)+\PADICVALUATION\left(\RHOOO_l\right)-\PADICVALUATION\left((\RHOOO p+l+1)_l\right)-l
\\
& \null\qquad\textstyle \GEQSLANT \gamma-1-\gamma-\frac{l-1}{p-1}+l=\frac{(l-1)(p-2)}{p-1} \GEQSLANT \frac{(p-2)^2}{(p-1)} > 0.
}
\end{proof}

\end{document}